\documentclass[conference]{IEEEtran}
\IEEEoverridecommandlockouts
\usepackage{cite}
\usepackage{amsmath,amssymb,amsfonts}
\allowdisplaybreaks
\usepackage{algorithmic}
\usepackage{graphicx}
\usepackage{subcaption}

\usepackage[utf8]{inputenc}
\usepackage[ruled, vlined, linesnumbered]{algorithm2e}
\SetKwInput{Input}{Input}
\SetKwInput{Output}{Output}
\usepackage{textcomp}
\usepackage{eurosym}
\usepackage{newfloat}
\DeclareFloatingEnvironment[placement={!ht},name=]{mylist}
\usepackage{float}
\usepackage{xcolor}
\usepackage{comment}
\usepackage{mathtools}
\usepackage{pgfplots}
\pgfplotsset{compat=1.18}
\usepgfplotslibrary{statistics,groupplots}

\def\BibTeX{{\rm B\kern-.05em{\sc i\kern-.025em b}\kern-.08em
    T\kern-.1667em\lower.7ex\hbox{E}\kern-.125emX}}

\begin{document}

\title{Machine Learning for Exact Time Series Aggregation in Generation Expansion Planning with Energy Storage\\ \thanks{Funded by the European Union (ERC, NetZero-Opt, 101116212). Views and opinions expressed are however those of the authors only and do not necessarily reflect those of the European Union or the European Research Council. Neither the European Union nor the granting authority can be held responsible for them.}
}

\author{\IEEEauthorblockN{Jakub Rybka, Luca Santosuosso, Thomas Klatzer and Sonja Wogrin}
\IEEEauthorblockA{\textit{Institute of Electricity Economics and Energy Innovation} \\
\textit{Research Center ENERGETIC} \\
\textit{Graz University of Technology} \\
Graz, Austria \\
\{j.rybka, luca.santosuosso, thomas.klatzer, wogrin\}@tugraz.at}
}
\maketitle

\newcommand\copyrighttext{%
  \footnotesize \textcopyright \the\year{} IEEE. Personal use of this material is permitted. Permission from IEEE must be obtained for all other uses, including reprinting/republishing this material for advertising or promotional purposes, collecting new collected works for resale or redistribution to servers or lists, or reuse of any copyrighted component of this work in other works.}

\newcommand\copyrightnotice{%
\begin{tikzpicture}[remember picture,overlay]
\node[anchor=south,yshift=10pt] at (current page.south) {\fbox{\parbox{\dimexpr0.75\textwidth-\fboxsep-\fboxrule\relax}{\copyrighttext}}};
\end{tikzpicture}%
}

\newcommand\IEEEcopyrighttext{\textbf{979-8-3195-3554-2/26/\$31.00 ©2026 IEEE}}
\newcommand\IEEEcopyrightnotice{%
\begin{tikzpicture}[remember picture,overlay]
\node[anchor=south west, xshift=1.6cm, yshift=1.3cm] 
at (current page.south west) {%
 \IEEEcopyrighttext};
\end{tikzpicture}%
}

\copyrightnotice 
\vspace{-10pt}

\begin{abstract}

This paper investigates a generation expansion planning (GEP) problem encompassing renewable, thermal, and storage technologies while simultaneously optimizing market participation, operational expenditures, and capital investment. 
To alleviate the computational burden of the GEP model, we propose a novel iterative time series aggregation (TSA) method that constructs a temporally aggregated counterpart of the original full-scale GEP model. 
Unlike traditional TSA methods, which are purely heuristic, our method enables the assessment of the optimality gap between the aggregated and full-scale models. 
Moreover, by leveraging machine learning–based estimates of the GEP model marginal costs, the algorithm guides TSA to construct an aggregated model that preserves the active constraints of its full-scale counterpart, which has been shown to yield exact temporal aggregation. 
Numerical results show that incorporating estimated marginal costs as clustering features substantially improves the quality of temporal aggregation compared with traditional TSA methods that rely solely on input data analysis.
\end{abstract}


\section{Introduction}

The generation expansion planning (GEP) problem is a classical optimization problem in power systems, focused on determining the optimal mix and sizing of power generation units to satisfy projected energy demand \cite{Koltsaklis2018}.
Historically, GEP has been primarily addressed by centralized system planners \cite{Lopez2020}.
However, in the context of liberalized electricity markets and increasing decentralization, this problem has gained renewed relevance for generation companies and emerging actors such as virtual power plants \cite{santosuosso2023economic},
whose objectives extend beyond cost minimization to the strategic design of an energy portfolio that maximizes profitability \cite{Sadeghi2017} while participating in multiple electricity markets \cite{santosuosso2024stochastic}.

The GEP problem is typically formulated over long-term planning horizons, often spanning multiple years \cite{Li2022-1}. 
Furthermore, the growing penetration of variable renewable energy (VRE) sources and flexible units, such as energy storage systems,
necessitates the explicit modeling of short-term operational dynamics to accurately capture the system evolution \cite{Levin2024}.
As a result, GEP often translates into large-scale, high-temporal-resolution optimization models that can become computationally demanding or even intractable \cite{Kotzur2021}.

To mitigate this computational challenge, time series aggregation (TSA) has emerged as a widely adopted approach \cite{Hoffmann2020}.
TSA reduces the full-scale GEP model to an aggregated model defined over a reduced set of representative time steps or clusters \cite{Pineda2018}.
When the number of representative steps is significantly smaller than the original time steps in the full-scale model, 
solving the aggregated model yields a significant computational advantage over its full-scale counterpart.

TSA methods are generally classified as \textit{a priori} or \textit{a posteriori} \cite{Hoffmann2020}.
Traditional a priori methods perform TSA based solely on the statistical features of the input time series to the GEP model,
typically employing standard clustering techniques such as K-Means \cite{Kotzur2018}, K-Medoids \cite{Teichgraeber2019}, or hierarchical clustering \cite{Moradi2023}. 
However, accurately representing the GEP input space does not guarantee an accurate representation of its output space,
meaning that the aggregated model output may still deviate significantly from that of the full-scale model \cite{Teichgraeber2022}.

Conversely, a posteriori TSA leverages structural information from the full-scale model to guide temporal aggregation.
For example, historical optimization runs can be used to estimate optimal decision variable values, which are then employed as clustering features \cite{Zhang2023}. 
Alternatively, daily subproblems may be solved to approximate the distribution of investment costs, which subsequently serve as features in K-Medoids \cite{Sun2019} or K-Means \cite{Li2022-2} clustering, rather than relying solely on input data analysis. 
A related approach extracts information from an initial a priori aggregated model to identify and emphasize extreme periods in the final aggregation \cite{Hilbers2023}.

Notably, previous research has demonstrated that a posteriori TSA methods capable of constructing an aggregated model that preserves the active constraints of the full-scale model achieve \textbf{exact} temporal aggregation \cite{Wogrin2023}; 
that is, the aggregated model yields the same optimal objective function value as the original full-scale model.
This theoretical result has been further extended to optimization models incorporating storage intertemporal constraints,
where TSA is further challenged by the need to maintain temporal chronology \cite{Klatzer2025}.
In practical applications, however, it is unrealistic to assume prior knowledge of the GEP model’s active constraints. Consequently, a critical research gap persists in the development of a practical a posteriori TSA method capable of constructing an aggregated model that approximates the active constraints of its full-scale counterpart,
targeting exact temporal aggregation while preserving chronology under intertemporal constraints.

This paper addresses this research gap through the following key contributions:
\begin{itemize}
    \item We formulate a GEP model that accounts for both capital investment and operational performance in the energy market while incorporating storage constraints. Building on this formulation, we propose a machine learning (ML) framework to estimate the marginal costs of the GEP model, i.e., the dual variables associated with the energy balancing constraints. These estimates are then used as clustering features to construct an aggregated model that preserves the active constraints of its full-scale counterpart, thereby explicitly targeting exact TSA.
    
    \item The proposed ML framework, which yields marginal cost estimates, is embedded within a novel \textit{a posteriori} TSA-based solution algorithm that provides explicit bounds on the approximation error introduced by the temporal aggregation procedure, thereby offering a clear and practical performance guarantee to the decision-maker.
    
    \item We benchmark the proposed TSA-based solution algorithm against state-of-the-art methods, including both traditional a priori TSA and recent a posteriori alternatives, to assess its potential for reducing computational complexity while yielding improved accuracy.
\end{itemize}

Finally, the performance of the proposed algorithm is validated using real-world Austrian hourly data on VRE capacity factors, energy demand, and energy prices for the year 2024.

The remainder of the paper is structured as follows.
Section \ref{s:method} introduces the proposed methodology.
Section \ref{s:res} presents the numerical results,
and Section \ref{s:con} concludes the paper.

\section{Methodology}\label{s:method}

This section presents the proposed methodology. Subsection \ref{sub:Agg_model} describes the aggregated GEP model. The ML framework and solution algorithm are detailed in Subsections \ref{sub:ML} and \ref{sub:Algorithm}, respectively.

\subsection{Aggregated Model}\label{sub:Agg_model}

Let $\mathcal{T}$, indexed by $t \in \{1, \dots, T\}$, denote the temporal horizon of the GEP problem.
A traditional full-scale formulation is provided in Appendix \ref{sub:Full_scale}.
As $\mathcal{T}$ grows, the full-scale model may become computationally intractable.
To address this, we introduce an aggregated GEP model defined over a reduced set of representative time steps $\mathcal{R}$, indexed by $r \in \{1, \dots, R\}$.

Let $\mathcal{T}_r \subseteq \mathcal{T}$ denote the set of consecutive time steps in $\mathcal{T}$ associated with representative step $r$, with cardinality $W_r \coloneqq |\mathcal{T}_r|$.
Moreover, let $\mathcal{G}$ and $\mathcal{S}$ denote the sets of generators, and storage units, respectively, indexed by $g$ and $s$.
All sets are indexed starting from $1$.

The charging and discharging efficiencies of the $s$-th storage unit are denoted by $\eta_s^\mathrm{c}$ and $\eta_s^\mathrm{d}$, respectively.
We denote by $\overline{E}_s$, $\overline{P}^\mathrm{d}_s$, $\overline{P}^\mathrm{c}_s$, and $\underline{E}_s$ the maximum state of charge (MWh), discharging capacity (MW), charging capacity (MW), and minimum state of charge (MWh) of the $s$-th storage unit, respectively.
The investment and operational costs of $g$ are denoted by $C^\mathrm{inv}_g$ (\euro/MW) and $C_g$ (\euro/MWh), respectively.
The penalty for non-supplied energy and the cost of discharging are $C^\mathrm{ns}$ (\euro/MWh) and $C^\mathrm{d}_s$ (\euro/MWh), respectively.
We denote by $\overline{B}$ the investment budget (\euro).
Finally, $\Delta$ denotes the time resolution (h) of the GEP problem.

Let $\hat{x}_g$ denote the capacity (MW) of generator $g$, $\hat{p}_{g,r}$ its power output (MW), and $\hat{e}^\mathrm{ns}_r$ the non-supplied energy (MWh) at the representative time step $r$.
We denote by $\hat o_r$ the energy (MWh) sold to the market at the representative time step $r$.
Moreover, $\hat e_{s,r}$ denotes the state of charge (MWh) of the storage unit $s$ at the representative time step $r$,
while $\hat p^\mathrm{c}_{s,r}$ and $\hat p^\mathrm{d}_{s,r}$ (MW) denote its charging and discharging power, respectively.
We introduce an auxiliary variable $\hat e_{s,0}$ to represent the initial state of charge of the storage unit $s$.

For all $r \in \mathcal{R}$, the aggregated input time series, denoted by $\hat{F}_{g,r}$, $\hat{D}_r$, and $\hat{\pi}_r$, are given by:
\begin{equation*}
    \hat{F}_{g,r} \coloneqq \sum_{t \in \mathcal{T}_r} \frac{F_{g,t}}{W_r}, \;\;\;\; \hat{D}_{r} \coloneqq \sum_{t \in \mathcal{T}_r} \frac{D_{t}}{W_r}, \;\;\;\; \hat{\pi}_{r} \coloneqq \sum_{t \in \mathcal{T}_r} \frac{\pi_t}{W_r},
\end{equation*}
where $F_{g,t}$, $D_t$, and $\pi_t$ denote capacity factors (p.u.), demand (MWh), and market signal time series (\euro/MWh), respectively.

We group the decision variables of the aggregated model into the set $\boldsymbol{\hat{z}}$ as follows:
\begin{equation} \label{eqn:decset_agg} 
\boldsymbol{\hat{z}} \coloneqq\left\{ \hat{x}_g, \hat{p}_{g,r}, \hat{p}^\mathrm{d}_{s,r}, \hat{p}^\mathrm{c}_{s,r}, \hat{e}^\mathrm{ns}_r, \hat{e}_{s,r}, \hat{o}_r\right\}_{g\in\mathcal{G}, r\in\mathcal{R}, s\in\mathcal{S}}, 
\end{equation}
while the aggregated objective function is defined as
\begin{align}
\label{eqn:obj_agg}
\!\!\!\!\!\!\hat{f}\!\left(\boldsymbol{\hat{z}}\right) \!\coloneqq\! \sum_{g \in \mathcal{G}} C_g^\mathrm{inv} \hat{x}_g \!+\! \sum_{r \in \mathcal{R}} \!W_r\! \bigg(& \!\Delta\! \sum_{g \in \mathcal{G}} C_g \hat{p}_{g,r} \!+\! \Delta \! \sum_{s \in \mathcal{S}} C_s^\mathrm{d} \hat{p}_{s,r}^\mathrm{d} 
\nonumber \\ &  + C^\mathrm{ns} \hat{e}_r^\mathrm{ns} - \hat\pi_r \hat{o}_r \bigg).
\end{align}

The \textbf{aggregated GEP model} is formulated as follows:
\begin{subequations} \label{eqn:model_agg} 
\begin{align} 
\!\!\!\!\min_{\boldsymbol{\hat{z}} \geq \boldsymbol{0}} \; & \hat{f}\left(\boldsymbol{\hat{z}}\right) \label{eqn:obj__agg} \\ 
\mathrm{s.t.} \; & \sum_{g \in \mathcal{G}} \hat{p}_{g,r} \Delta \!+\! \sum_{s \in \mathcal{S}} \left(\hat{p}_{s,r}^\mathrm{d} \!-\! \hat{p}_{s,r}^\mathrm{c}\right) \!\Delta \!+\! \hat{e}_r^\mathrm{ns} \!-\! \hat{o}_r = \hat{D}_r,\!\! \forall r, \label{eqn:bal_agg} \\ 
& \hat{p}_{g,r} \leq \hat{F}_{g,r} \, \hat{x}_g, \; \forall g, \forall r, \label{eqn:gen_agg} \\
& \hat{e}_{s,r} = \hat{e}_{s,r-1} + \left(\eta_s^\mathrm{c} \hat{p}_{s,r}^\mathrm{c} - \frac{\hat{p}_{s,r}^\mathrm{d}}{\eta_s^\mathrm{d}} \right) \Delta \, W_r, \; \forall s,  \forall r, \label{eqn:soc_dyn_agg} \\ 
& \hat{e}_{s,0} = \hat{e}_{s,R} = \underline{E}_s, \; \forall s,\label{eqn:soc_agg_end} \\
& \underline{E}_s \leq \hat{e}_{s,r} \leq \overline{E}_s, \; \forall s, \forall r, \label{eqn:soc_lim_agg} \\ 
& \hat{p}_{s,r}^\mathrm{c} \leq \overline{P}_s^\mathrm{c}, \; \forall s, \forall r, \label{eqn:storage_charge_agg}\\
& \hat{p}_{s,r}^\mathrm{d} \leq \overline{P}_s^\mathrm{d}, \; \forall s, \forall r, \label{eqn:storage_discharge_agg} \\ 
& \sum_{g \in \mathcal{G}} C_g^\mathrm{inv} \hat{x}_g \leq \overline{B}. \label{eqn:inv_lim_agg}
\end{align} 
\end{subequations}

In \eqref{eqn:model_agg}, objective \eqref{eqn:obj__agg} scales operational costs by $W_r$, while \eqref{eqn:bal_agg} ensures energy balance. 
Generation limits are enforced by \eqref{eqn:gen_agg}, and storage dynamics by \eqref{eqn:soc_dyn_agg}--\eqref{eqn:soc_lim_agg}. 
Finally, \eqref{eqn:storage_charge_agg}, \eqref{eqn:storage_discharge_agg}, and \eqref{eqn:inv_lim_agg} bound charging, discharging, and cumulative investment.
Notably, when $\mathcal{R}=\mathcal{T}$, the model is equivalent to the full-scale formulation provided in Appendix \ref{sub:Full_scale}.

\subsection{Machine Learning for Marginal Cost Estimation}\label{sub:ML}

Klatzer et al. \cite{Klatzer2025} demonstrate that aggregating consecutive time steps with identical active constraints yields \textbf{exact} temporal aggregation. 
Because these constraints are unknown a priori, we propose a \textbf{marginal-cost-based} (\textbf{MCB}) \textbf{TSA method}. 
This method employs an ML framework to estimate full-scale marginal costs as clustering features, ensuring the aggregated model \eqref{eqn:model_agg} preserves the active constraints of the full-scale model \eqref{eqn:model}.

The proposed ML framework takes the GEP input time series $\Omega^\mathcal{T}$, comprising capacity factors, demand, and energy prices, and estimates its unknown true distribution $P\left(\Omega^\mathcal{T}\right)$ via Kernel Density Estimation \cite{Chen2017}, yielding the estimate $\bar{P}\left(\Omega^\mathcal{T}\right)$.
An artificial time series $\Omega^\mathcal{K}$ is generated by drawing $K$ samples from $\bar{P}\left(\Omega^\mathcal{T}\right)$
and used to solve a reduced version of the full-scale GEP model \eqref{eqn:model}, defined over $\mathcal{K}$ instead of $\mathcal{T}$. 
The outputs of this reduced GEP model are then used to train a Random Forest classifier \cite{Parmar2018}, which is subsequently employed to generate estimates of the marginal costs $\boldsymbol{\bar \mu}$ over the full planning horizon $\mathcal{T}$.
The classifier is configured with $M = 100$ trees and a maximum depth of $20$. The bootstrap aggregation mitigates overfitting, while the chosen depth enables the capture of marginal cost variations without sacrificing generalization. With a training complexity of $O(M \cdot K \log K)$, the computational overhead of the ML framework remains negligible compared to that of the full-scale GEP model.

The proposed ML framework is summarized as follows:
\begin{enumerate}\label{ML_workflow}
    \item \textbf{Input:} $\Omega^\mathcal{T} \coloneqq \left\{F_{g,t}, D_t, \pi_t \right\}_{g\in\mathcal{G}, t\in\mathcal{T}}$.
    \item \textbf{Density estimation:} Estimate the input data distribution $\bar P\left(\Omega^\mathcal{T}\right) \approx P\left(\Omega^\mathcal{T}\right)$ via Kernel Density Estimation \cite{Chen2017}.
    \item \textbf{Sampling:} Draw $K$ samples from $\bar P\left(\Omega^\mathcal{T}\right)$ to generate a time series
    $\Omega^\mathcal{K} \coloneqq \left\{\bar P(\Omega^\mathcal{T})_k\right\}_{k \in \mathcal{K}}$, with $K \ll |\mathcal{T}|$.
    \item \textbf{Marginal costs extraction:} Solve the GEP model \eqref{eqn:model} for $\Omega^\mathcal{K}$ to obtain the marginal cost estimates $\boldsymbol{\tilde \mu} \coloneqq \left\{\tilde\mu_k\right\}_{k \in \mathcal{K}}$.
    \item \textbf{Training:} Train a Random Forest \cite{Parmar2018} classifier $\mathcal{C}: \Omega^{\mathcal{T}} \to \boldsymbol{\bar \mu}$ using the labeled pairs $\left(\Omega^\mathcal{K}, \boldsymbol{\tilde \mu}\right)$.
    \item \textbf{Inference:} Predict the marginal costs $\boldsymbol{\bar \mu} \coloneqq \left\{\bar\mu_t\right\}_{t \in \mathcal{T}}$ for the full planning horizon $\mathcal{T}$ using $\mathcal{C}$.
\end{enumerate}

\subsection{Marginal-Cost-Based Time Series Aggregation}\label{sub:Algorithm}

 The aggregated model \eqref{eqn:model_agg} solution provides a lower bound ($f^{\text{LB}}$) to the optimal objective function value ($f^\star$) of the full-scale model \eqref{eqn:model} \cite{Santosuosso2025}. 
Projecting the aggregated model solution onto the feasible space of \eqref{eqn:model} then yields an upper bound ($f^{\text{UB}}$).
In Algorithm \ref{alg:ts_aggregation}, this bounding procedure is combined with the ML framework to iteratively refine the aggregated model \eqref{eqn:model_agg} toward exact TSA, while monitoring the achieved optimality gap at each iteration.

\begin{algorithm}[ht]
    \caption{Time Series Aggregation Guided by Marginal Cost Estimates}
    \label{alg:ts_aggregation}
    
    \Input{$\Omega^\mathcal{T}$, $R^{0}$, $K$, $N^{\mathrm{top}}$, $\epsilon^{\mathrm{target}}$, $\overline{R}$, $\delta_R$.}
    \Output{Objective function bounds $f^{\text{UB}}$ and $f^{\text{LB}}$.}
    
    Predict $\boldsymbol{\bar \mu}$ using the ML framework of Subsection \ref{ML_workflow}\;
    
    Construct $\mathcal{T}^{\mathrm{top}}$ consisting of the $N^{\mathrm{top}}$ time steps in $\mathcal{T}$ with the highest net demand $\tilde{D}_t$, calculated as in \eqref{eqn:net_dem}\;

    $\mathcal{T}^-\coloneqq \mathcal{T}\setminus\mathcal{T}^{\mathrm{top}}$\;

    $\left\{\mathcal{R}^{\mathrm{top}}, \mathcal{T}_{\mathcal{R}^{\mathrm{top}}}\right\} \gets$ For each $\tilde{t} \in \mathcal{T}^{\mathrm{top}}$, create a representative time step $\tilde{r}$ with $\mathcal{T}_{\tilde{r}} \coloneqq \left\{\tilde{t}\right\}$, then include $\tilde{r}$ in $\mathcal{R}^{\mathrm{top}}$ and $\mathcal{T}_{\tilde{r}}$ in $\mathcal{T}_{\mathcal{R}^{\mathrm{top}}}$.
    
    \textbf{Initialization:} $\epsilon \gets+ \infty, \quad R \gets R^{0}$\;
    
    \While{$\epsilon > \epsilon^{\mathrm{target}}$ \textbf{and} $R < \overline{R}$}{
    
        $\left\{\mathcal{R}^-,\mathcal{T}_{\mathcal{R}^-}\right\} \gets$ Perform chronological hierarchical clustering \cite{Pineda2018}, using $\left\{\bar\mu_t\right\}_{t \in \mathcal{T}^-}$ as clustering features, $\left|\bar\mu_t-\bar\mu_{t'}\right|$ as distance metric between consecutive time steps $t$ and $t'$, and $R$ as the desired number of representative time steps\;
        $\mathcal{R}\gets \mathcal{R}^- \cup \mathcal{R}^{\mathrm{top}}$\;
        $\mathcal{T}_{\mathcal{R}}\gets\mathcal{T}_{\mathcal{R}^-} \cup \mathcal{T}_{\mathcal{R}^{\mathrm{top}}}$\;
        
        $\boldsymbol{\hat{z}}^\star \gets$ Solve the aggregated GEP model \eqref{eqn:model_agg} over $\mathcal{R}\space$\;
        
        $f^{\text{LB}} \gets \hat{f}\left(\boldsymbol{\hat{z}}^\star\right)$ \tcp*{Lower bound from the aggregated objective function \eqref{eqn:obj_agg}}

        $\boldsymbol{\bar{z}}^\star \gets$ Solve \eqref{eqn:model} with generation capacity variables fixed to their values in $\boldsymbol{\hat{z}}^\star$\;
        
        $f^{\text{UB}} \gets f(\boldsymbol{\bar{z}}^\star)$ \tcp*{Upper bound}
        
        $\epsilon \gets 100 \, \frac{\left(f^{\text{UB}} - f^{\text{LB}}\right)}{f^{\text{UB}}}$\tcp*{Optimality gap}
        $R \gets R + \delta_R$\;
    }
\end{algorithm}

The algorithm uses the marginal cost estimates $\boldsymbol{\bar{\mu}}$ as clustering features for TSA, within the chronological hierarchical clustering algorithm of \cite{Pineda2018}.
It terminates when the optimality gap $\epsilon$, defined as the relative difference between the upper and lower bounds, falls below the target $\epsilon^\mathrm{target}$ or when the maximum aggregated model size $\overline{R}$ is reached.
At each iteration, the number of representative steps $R$ is increased from the initial value $R^0$ by $\delta_R$, progressively refining the objective function bounds.
To improve aggregation fidelity, full resolution is maintained for a protected set $\mathcal{T}^{\mathrm{top}}$ comprising the $N^{\mathrm{top}}$ time steps with the highest net demand.
The net demand $\tilde{D}_t$ is estimated from the generation capacities $\tilde{x}_g$ obtained via the ML framework in Subsection \ref{ML_workflow}, as follows:
\begin{equation}\label{eqn:net_dem}
    \tilde{D}_t  = D_t - \Delta\sum_{g\in\mathcal{G}} \tilde x_g \, F_{g,t}, \; \forall t.
\end{equation}

The proposed Algorithm \ref{alg:ts_aggregation} exhibits two notable features.
First, it provides a measure of the optimality gap achieved by the TSA method at each iteration \cite{Santosuosso2025}, thereby offering a formal performance guarantee to the decision-maker, whereas traditional a priori TSA methods are purely heuristic.
Second, it leverages structural information from the full-scale GEP model, specifically the marginal cost estimates, to guide TSA toward exact temporal aggregation, whereas conventional a priori TSA methods rely solely on input data analysis.

\section{Numerical Results}\label{s:res}

This section presents a performance evaluation of the proposed MCB TSA in comparison with several state-of-the-art TSA methods.
In Subsection~\ref{sub:CM}, we consider a cost minimization problem in which market participation is disregarded in the GEP model,
while Subsection~\ref{sub:MP} considers a GEP problem that includes market participation.

The case study relies on historical open-source data on demand, capacity factors, and energy prices in Austria for 2024, obtained from the ENTSO-E Transparency Platform \cite{Hirth2018}.
These data are perturbed across 10 distinct scenarios using uniform noise drawn from the interval $[0.8, 1.2]$ to assess the sensitivity of the analyzed TSA methods to variations in the input data.
The sensitivity across scenarios is reported using boxplots, where the boxes denote the interquartile range, the median is shown in blue, and the whiskers encompass all data points.
Each case study includes a single storage unit with varying energy-to-power ratios (in hours), $\overline{E}_s / \overline{P}^\mathrm{d}_s$, along with one renewable and one thermal generation unit.

We assume generation costs of $130$, $2.5$, and $1$ (\euro/MWh) for thermal, wind, and PV units, respectively.
The investment costs are $10^5$ (\euro/MW) for thermal units and $8 \times 10^4$ (\euro/MW) for VRE.
The non-supplied energy cost is $5 \times 10^3$ (\euro/MWh), and the storage discharging cost is $1.5$ (\euro/MWh).
For the storage unit, we set $\eta_s^\mathrm{c} = 0.9$, $\eta_s^\mathrm{d} = 0.9$, $\overline{P}^\mathrm{d}_s =200$ (MW), $\overline{P}^\mathrm{c}_s = 200$ (MW), $\underline{E}_s = 0$ (MWh).
The investment budget is $3 \times 10^8$ (\euro).

All GEP models are defined over 8,736 time steps at an hourly resolution.
In Algorithm~\ref{alg:ts_aggregation}, we set $\epsilon^{\mathrm{target}} = 1\%$, $\overline{R} = 900$, $R^{0} = 400$, $\delta_R = 100$, $K = 500$, and $N^{\mathrm{top}} = 100$.

To balance accuracy and computational complexity, the optimal values of $K$ and $N^{\mathrm{top}}$ were determined through a sensitivity analysis.
Results indicate that $K < 500$ fails to adequately capture the input time series distribution, leading to inaccurate marginal cost representations, whereas $K > 500$ yields negligible improvement.
Similarly, $N^{\mathrm{top}} < 100$ underestimates investment decisions by omitting extreme periods, while larger values provide no significant additional benefit.

The key differences between the proposed MCB TSA method and the benchmark methods are summarized in Table~\ref{tab:methods}.
For consistency, all methods are implemented within Algorithm \ref{alg:ts_aggregation}, with the adopted clustering features and techniques varying according to the specific TSA method.
Specifically, \cite{Pineda2018} and \cite{Moradi2023} propose a priori TSA methods that rely solely on the input data to identify representative hours \cite{Pineda2018} or days \cite{Moradi2023}.
We further compare our method with a posteriori TSA methods that exploit structural information from the GEP model through estimates of the optimal primal variable values, such as investment decisions \cite{Zhang2023} or net demand \cite{Hilbers2023},
whereas our method instead employs estimates of the dual information, namely the marginal costs, as clustering features.

\begin{table}[t]
\caption{Comparison of the analyzed TSA methods.}
\centering
\label{tab:methods}
\renewcommand{\arraystretch}{1.2}
\begin{tabular}{|l|c|c|c|}
\hline
\textbf{Method} & \textbf{Clustering features} & \textbf{Clustering}  & \textbf{Type of TSA method} \\ \hline
\cite{Pineda2018} & Input data & CHC &  A priori \\ \hline
\cite{Moradi2023} & Input data & HC  & A priori \\ \hline
\cite{Zhang2023} & Primal variables & K-Medoids  & A posteriori \\ \hline
\cite{Hilbers2023} & Primal variables & K-Medoids & A posteriori \\ \hline
\textbf{MCB} & \textbf{Dual variables} & \textbf{CHC} & \textbf{A posteriori} \\ \hline
\multicolumn{4}{l}{\scriptsize \textbf{Legend}: Chronological hierarchical clustering (CHC), hierarchical clustering (HC).} 
\end{tabular}
\end{table}

The GEP models were solved on an Intel Xeon W5 workstation with 256~GB of RAM using Gurobi 13.0.0.

\subsection{Cost Minimization Problem}\label{sub:CM}

\begin{figure}[ht]

    \centering
    \includegraphics[width=1\linewidth]{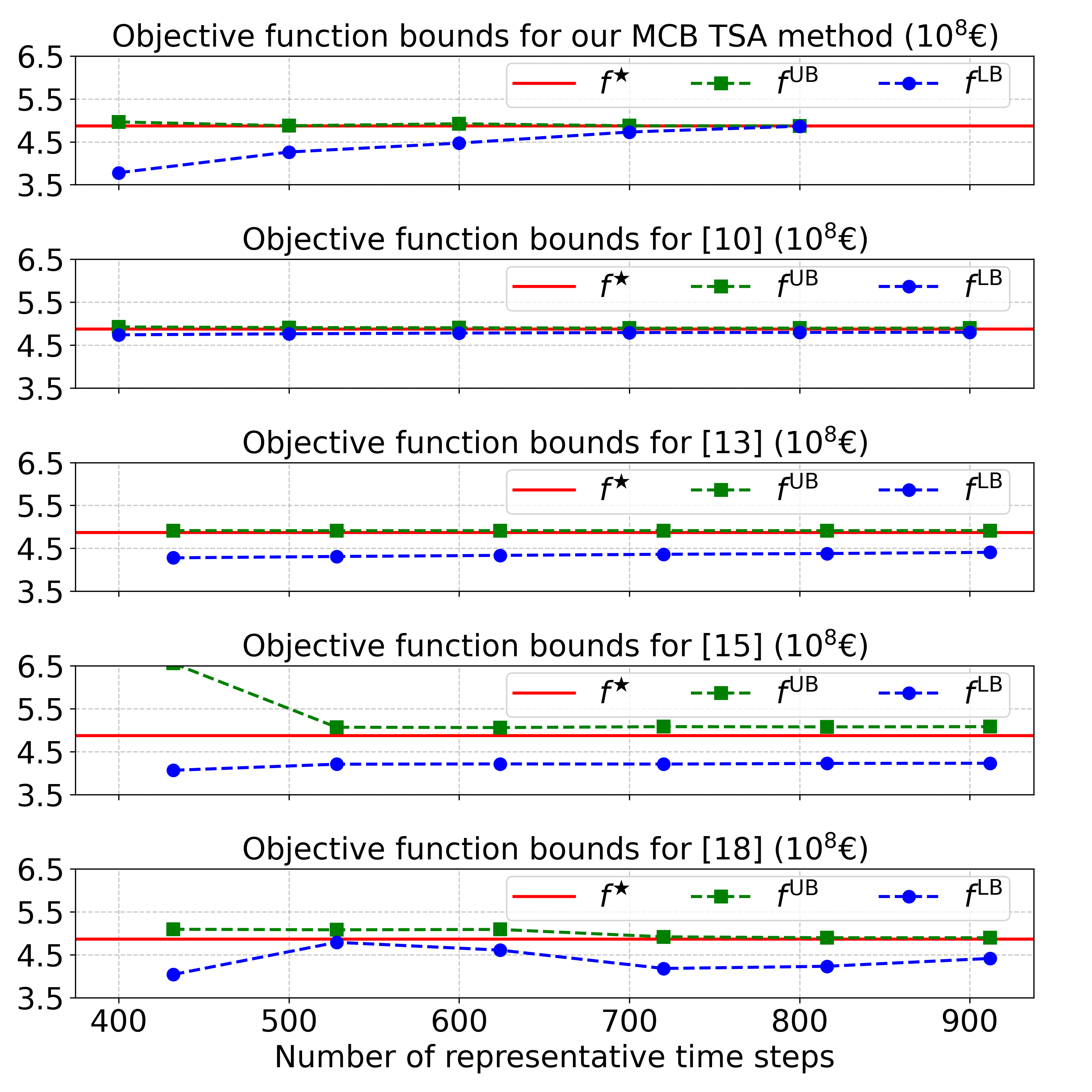}
    \captionsetup{skip=-3.5pt}
    \caption{Objective function bounds for a stylized cost minimization problem with a wind unit as the only VRE source and a storage unit with a 4h energy-to-power ratio.}
    \label{fig:conv_cm}
\end{figure}

The observed objective function bounds are reported in Fig.~\ref{fig:conv_cm}. Our MCB TSA method outperforms all benchmark methods, converging to the desired optimality gap with $800$ representative steps.
Notably, the TSA method of \cite{Pineda2018} yields the smallest initial optimality gap ($3.7\%$), observed at the first iteration of Algorithm \ref{alg:ts_aggregation}, but fails to reach the target gap of $1\%$.
Since our method differs from that of \cite{Pineda2018} in the clustering features employed (estimated marginal costs versus input time series), this result further underscores the effectiveness of marginal costs as proxies for exact TSA.

\begin{figure}[t]
     \begin{subfigure}{\linewidth}
\hspace{-0.55cm}
  \centering
  \begin{tikzpicture}
    \begin{groupplot}[
      group style={
        group size=4 by 1,
        horizontal sep=0.7cm,
        x descriptions at=edge bottom,
      },
      width=0.355\linewidth,
      height=3.9cm,
      boxplot/draw direction=y,
      xtick={1,...,5},
      xticklabels={MCB,\cite{Pineda2018},\cite{Moradi2023},\cite{Zhang2023},\cite{Hilbers2023}},
      xticklabel style={rotate=90, anchor=east, font=\scriptsize},
      ymajorgrids,
      grid style={dashed, gray!30},
      boxplot/every median/.style={thick,draw=blue},      boxplot/whisker range=5000,
      title style={font=\scriptsize\bfseries, yshift=-3pt},
      tick label style={font=\scriptsize},
      boxplot/every box/.style={fill=gray!30},
      boxplot/draw/marker/.style={draw=none},
    ]
      \nextgroupplot[title={\!\!\!Opt. gap (\%)},title style={yshift=-4pt}, ymin=-10.00, ymax=36.00, ytick distance=5.00]
      \draw[red, thick] (axis cs:0.5,0) -- (axis cs:5.5,0);
        \addplot+[black, boxplot] coordinates {(1,0.1999533909250805) (1,0.0963310559602164) (1,0.2009658446279618) (1,0.1585464597702612) (1,0.153251021649218) (1,0.2404219575997979) (1,0.0890992372471653) (1,0.1888606669976263) (1,0.0717776851460811) (1,0.1799246514860648)};
        \addplot+[black, boxplot] coordinates {(2,1.6586506712544955) (2,1.1958109920478817) (2,1.5177076561169722) (2,1.1939324388626489) (2,1.236466385577864) (2,1.148834295643233) (2,1.2684221776564624) (2,1.1263257890965652) (2,1.6702992536659425) (2,1.671504080665127)};
        \addplot+[black, boxplot] coordinates {(3,4.697369281986586) (3,5.528023230684044) (3,6.253498010021986) (3,5.774776259811776) (3,6.251665068948155) (3,5.974523679129873) (3,5.365412376139924) (3,5.459291896270661) (3,5.36681544653204) (3,5.156430447616984)};
        \addplot+[black, boxplot] coordinates {(4,4.67699274654502) (4,2.5144729338590723) (4,4.035225826204951) (4,2.8174944974794687) (4,5.323772619491581) (4,5.210022992520874) (4,2.873184826984238) (4,4.055930711145842) (4,2.755911140761284) (4,4.62145887231074)};
        \addplot+[black, boxplot] coordinates {(5,8.482642306302195) (5,5.036292761259263) (5,3.8762701850731) (5,3.3463551206788087) (5,3.812650203998944) (5,2.705603333826047) (5,4.110002606857977) (5,2.5502411652659758) (5,3.8564549026661816) (5,0.2578000592637784)};
        \node[draw,fill=white,font=\scriptsize,anchor=north west,inner xsep=2pt,inner ysep=1pt] at (rel axis cs:0.03,0.97){\raisebox{0.5ex}{\textcolor{red}{\rule{0.2cm}{0.7pt}}}~F-S};
      \nextgroupplot[title={\!\!\!VRE inv. error\!\! (\%)},title style={yshift=-4pt}, ymin=-10.00, ymax=36.00, ytick distance=5.00]
      \draw[red, thick] (axis cs:0.5,0) -- (axis cs:5.5,0);
        \addplot+[black, boxplot] coordinates {(1,0.7771530274705073) (1,-1.3521266994526209) (1,-1.394254942368242) (1,-0.7537812908294046) (1,1.3798480364813577) (1,-0.2041699237273174) (1,-2.097819103284261) (1,-1.587212527523245) (1,0.7086717517534656) (1,0.3964144869443343)};
        \addplot+[black, boxplot] coordinates {(2,5.697727229476463) (2,6.949849421449541) (2,9.4709908586898) (2,5.846894192607204) (2,5.010033134904069) (2,4.446023603690911) (2,4.879038840571958) (2,3.1478949711709454) (2,7.393016633904096) (2,8.831513659177002)};
        \addplot+[black, boxplot] coordinates {(3,20.423505226649336) (3,25.81951343706345) (3,33.29021386337799) (3,24.41833993558788) (3,29.2077601308722) (3,28.10210260436991) (3,22.7034014098001) (3,26.851254294761663) (3,23.266887616424334) (3,23.66846756215855)};
        \addplot+[black, boxplot] coordinates {(4,2.6625409658857677) (4,0.1577953161660023) (4,2.851405905195555) (4,5.177679105311126) (4,7.783306272905113) (4,1.5518264427166042) (4,4.979745450062624) (4,4.471198644029029) (4,4.832038057567205) (4,2.198835354178948)};
        \addplot+[black, boxplot] coordinates {(5,-6.292436756230131) (5,-6.579245178229785) (5,6.8300115976584665) (5,-1.0014837302454644) (5,3.200781543709984) (5,-4.160809388932612) (5,1.9266329500327983) (5,-2.4284813412762607) (5,-6.051476740768253) (5,6.414151448543426)};
        \node[draw,fill=white,font=\scriptsize,anchor=north west,inner xsep=2pt,inner ysep=1pt] at (rel axis cs:0.03,0.62){\raisebox{0.5ex}{\textcolor{red}{\rule{0.2cm}{0.7pt}}}~F-S};
      \nextgroupplot[title={$\;$ Th. inv. error\! (\%)}, title style={yshift=-4pt},ymin=-10.00, ymax=36.00, ytick distance=5.00]
      \draw[red, thick] (axis cs:0.5,0) -- (axis cs:5.5,0);
        \addplot+[black, boxplot] coordinates {(1,-2.8440992153679194) (1,-0.1605526409992001) (1,-2.4144871705821456) (1,-1.5529303010936073) (1,-1.351525700793489) (1,-0.8243023929435442) (1,-0.3899750737699649) (1,-1.0836334129624468) (1,-0.5752476751718779) (1,-1.529959847055247)};
        \addplot+[black, boxplot] coordinates {(2,-4.414954694559979) (2,-1.710831727336684) (2,-5.3161995500579176) (2,-4.19395893892631) (2,-2.5902085489186035) (2,-1.3390461968559462) (2,-3.321424949233556) (2,-2.62182040795836) (2,-4.777914141584624) (2,-5.785147876032964)};
        \addplot+[black, boxplot] coordinates {(3,-2.813133787546712) (3,-0.9718529575275978) (3,-2.4144871705821456) (3,-2.069716951867504) (3,-2.189289979986002) (3,-1.4628151162716585) (3,-1.883632821499956) (3,-1.6915119910122765) (3,-0.8906389702574241) (3,-1.98990453407458)};
        \addplot+[black, boxplot] coordinates {(4,-8.619673612179787) (4,-6.214444891574016) (4,-8.777168315278837) (4,-7.32316126990961) (4,-8.269620296829679) (4,-9.748652586345234) (4,-6.428800295387332) (4,-7.356134709025635) (4,-6.392041997361384) (4,-9.195801096279418)};
        \addplot+[black, boxplot] coordinates {(5,0.0457653202726832) (5,-1.9128181004990044) (5,-2.994443817752481) (5,0.0848865180424159) (5,-1.8786322456728517) (5,0.3737821453354463) (5,-0.1598389746471359) (5,-0.6172438176405367) (5,-0.7963574121914742) (5,-1.529959847055247)};
        \node[draw,fill=white,font=\scriptsize,anchor=north west,inner xsep=2pt,inner ysep=1pt] at (rel axis cs:0.03,0.97){\raisebox{0.5ex}{\textcolor{red}{\rule{0.2cm}{0.7pt}}}~F-S};
      \nextgroupplot[title={Work units (-)},title style={yshift=-4pt}, ymin=0.00, ymax=0.92, ytick distance=0.10]
        \draw[ red, thick] (axis cs:0.5,0.60) -- (axis cs:5.5,0.60);
        \addplot+[black, boxplot] coordinates {(1,0.0862018054148293) (1,0.0959958147896659) (1,0.0954816459722708) (1,0.0942292134202117) (1,0.0953680886076462) (1,0.092785120499355) (1,0.0938334266325325) (1,0.0891414227897463) (1,0.0953069520165901) (1,0.0940017986577605)};
        \addplot+[black, boxplot] coordinates {(2,0.0853370132528683) (2,0.0861348342808303) (2,0.0844833081445726) (2,0.0833930510874025) (2,0.0886945123450412) (2,0.0893428103959477) (2,0.0852362960468064) (2,0.0839238967742472) (2,0.0841239653358634) (2,0.0879211919735591)};
        \addplot+[black, boxplot] coordinates {(3,0.3490068695009812) (3,0.3853694658008281) (3,0.4336536864282341) (3,0.4417000860524128) (3,0.2821887513368087) (3,0.3919915157320399) (3,0.365537072534065) (3,0.3679188827407001) (3,0.3280192161883086) (3,0.304294793852522)};
        \addplot+[black, boxplot] coordinates {(4,0.5907644430796305) (4,0.41660769780476886) (4,0.39520470301310223) (4,0.8828331629435221) (4,0.4286106427510579) (4,0.41495569547017414) (4,0.42508212725321454) (4,0.4464809894561768) (4,0.4085397720336914) (4,0.4399905204772949)};
        \addplot+[black, boxplot] coordinates {(5,0.2974098033801982) (5,0.2806739677663285) (5,0.3108363777735084) (5,0.2633621795309777) (5,0.2853301826410427) (5,0.2632906909948038) (5,0.2575792203359809) (5,0.251481860326084) (5,0.3251780337712547) (5,0.2891010317921723)};
        \node[draw,fill=white,font=\scriptsize,anchor=north west,inner xsep=2pt,inner ysep=1pt] at (rel axis cs:0.03,0.97){\raisebox{0.5ex}{\textcolor{red}{\rule{0.2cm}{0.7pt}}}~F-S};
    \end{groupplot}
  \end{tikzpicture}
  \captionsetup{skip=-3.5pt}
  \caption{Case study with 4h energy-to-power ratio storage.}
  \label{fig:ess4_solar}
\end{subfigure}
     \begin{subfigure}{\linewidth}
\hspace{-0.55cm}
  \centering
  \begin{tikzpicture}
    \begin{groupplot}[
      group style={
        group size=4 by 1,
        horizontal sep=0.7cm,
        x descriptions at=edge bottom,
      },
      width=0.355\linewidth,
      height=3.9cm,
      boxplot/draw direction=y,
      xtick={1,...,5},
      xticklabels={MCB,\cite{Pineda2018},\cite{Moradi2023},\cite{Zhang2023},\cite{Hilbers2023}},
      xticklabel style={rotate=90, anchor=east, font=\scriptsize},
      ymajorgrids,
      grid style={dashed, gray!30},
      boxplot/every median/.style={thick,draw=blue},      boxplot/whisker range=5000,
      title style={font=\scriptsize\bfseries, yshift=-3pt},
      tick label style={font=\scriptsize},
      boxplot/every box/.style={fill=gray!30},
      boxplot/draw/marker/.style={draw=none},
    ]
    
      \nextgroupplot[title={Opt. gap (\%)},title style={yshift=-4pt}, ymin=-20.00, ymax=31.00, ytick distance=5.00]
      \draw[red, thick] (axis cs:0.5,0) -- (axis cs:5.5,0);
        \addplot+[black, boxplot] coordinates {(1,0.1170442165229036) (1,0.1727321105462991) (1,0.1272593853832061) (1,0.1973892069815241) (1,0.2241302977009261) (1,0.072645953091067) (1,0.1478814061589729) (1,0.0958235269259358) (1,0.1095576177610949) (1,0.1472131056601348)};
        \addplot+[black, boxplot] coordinates {(2,6.873138127696218) (2,6.282442997074974) (2,6.544730234942228) (2,8.207157173495954) (2,6.357143428969804) (2,6.991501234409289) (2,7.099859432380148) (2,7.008921249392987) (2,5.212998561462193) (2,7.037516746990422)};
        \addplot+[black, boxplot] coordinates {(3,3.6508082347679536) (3,4.822120446097651) (3,4.8387825684501165) (3,5.375184690272398) (3,4.343681008626828) (3,5.561667916066233) (3,6.096307923397604) (3,5.77292891373368) (3,4.922366741137086) (3,5.010859393323791)};
        \addplot+[black, boxplot] coordinates {(4,4.324803133348559) (4,4.455019695294933) (4,2.924927494557193) (4,3.980542363485331) (4,3.490880976758544) (4,5.078141307675399) (4,4.386526934097651) (4,3.9546560029392928) (4,4.2994824718341205) (4,3.946706947332304)};
        \addplot+[black, boxplot] coordinates {(5,0.4797861620484324) (5,2.2252918802698964) (5,6.365010343951419) (5,4.079364959335837) (5,3.0206694826777603) (5,5.883111327712931) (5,3.576827237269117) (5,4.1045271097639) (5,3.117135069081892) (5,1.9711786562899385)};
        \node[draw,fill=white,font=\scriptsize,anchor=north west,inner xsep=2pt,inner ysep=1pt] at (rel axis cs:0.03,0.97){\raisebox{0.5ex}{\textcolor{red}{\rule{0.2cm}{0.7pt}}}~F-S};
      \nextgroupplot[title={\!\!\!VRE inv. error\! (\%)}, title style={yshift=-4pt},ymin=-20.00, ymax=31.00, ytick distance=5.00]
      \draw[red, thick] (axis cs:0.5,0) -- (axis cs:5.5,0);
        \addplot+[black, boxplot] coordinates {(1,3.041701809299244) (1,-7.17931101023731) (1,2.948165318116245) (1,-6.27190314982273) (1,-7.402975061932593) (1,3.0958887711978647) (1,4.518051909871524) (1,-2.894359511073705) (1,3.8229701492247794) (1,5.132679433710347)};
        \addplot+[black, boxplot] coordinates {(2,3.3747160527658906) (2,3.3930364580398185) (2,4.478246314722744) (2,8.042549167545145) (2,5.639725534083928) (2,3.6582505469006272) (2,7.225997189998353) (2,3.610680755651362) (2,1.9278466600519144) (2,2.947819645969507)};
        \addplot+[black, boxplot] coordinates {(3,8.98476373006466) (3,17.494870975809455) (3,15.932663160389923) (3,17.433386497176212) (3,8.199835195355146) (3,18.997151980419044) (3,25.64995386695178) (3,22.236731722457005) (3,18.12911695686656) (3,19.977455708101544)};
        \addplot+[black, boxplot] coordinates {(4,5.162872635745577) (4,1.4663233053729) (4,2.0210076383333284) (4,2.9122462200749983) (4,1.7332582593936183) (4,5.636761179560508) (4,5.797036026849231) (4,4.273495964837421) (4,3.877643774832141) (4,-1.5165463378081387)};
        \addplot+[black, boxplot] coordinates {(5,1.2378854254075122) (5,-8.17734487137492) (5,3.1662672396637808) (5,2.425785385725339) (5,0.8190587434515596) (5,1.5884628134157044) (5,-0.2378394071414026) (5,-8.651292490497724) (5,-3.571581693378546) (5,6.548695417822414)};
        \draw[red, thick] (axis cs:0.5,0) -- (axis cs:5.5,0);
        \node[draw,fill=white,font=\scriptsize,anchor=north west,inner xsep=2pt,inner ysep=1pt] at (rel axis cs:0.03,0.15){\raisebox{0.5ex}{\textcolor{red}{\rule{0.2cm}{0.7pt}}}~F-S};
      \nextgroupplot[title={$\;$ Th. inv. error\! (\%)},title style={yshift=-4pt}, ymin=-20.00, ymax=31.00, ytick distance=5.00]
      \draw[red, thick] (axis cs:0.5,0) -- (axis cs:5.5,0);
        \addplot+[black, boxplot] coordinates {(1,-0.0725993257703875) (1,-0.0267147097045566) (1,-0.3945710307392734) (1,-0.0481949150173658) (1,-0.3109788749864032) (1,-0.9640290901604004) (1,-0.4465818082832914) (1,-1.193472474285297) (1,-0.0113151667103533) (1,-0.6681888638643384)};
        \addplot+[black, boxplot] coordinates {(2,-13.930970198862006) (2,-14.040513761183112) (2,-13.296371665140857) (2,-15.820315946659802) (2,-13.56770179541143) (2,-16.961888536899224) (2,-15.581345630753532) (2,-15.148254047654916) (2,-13.249086935854931) (2,-15.695905934725952)};
        \addplot+[black, boxplot] coordinates {(3,-1.3361174416796917) (3,-1.742292957989126) (3,-1.1674690301160704) (3,-2.1800252150659403) (3,-1.729116378000417) (3,-4.729427245802584) (3,-2.7591310686072235) (3,-2.3319004096276497) (3,-0.38039965617495) (3,-3.157591106808264)};
        \addplot+[black, boxplot] coordinates {(4,-11.360523208557945) (4,-11.733921216178093) (4,-7.610408249709866) (4,-9.82742445643316) (4,-9.270428505078208) (4,-13.527827972535944) (4,-11.90117171134381) (4,-10.80932591938817) (4,-13.002696944990385) (4,-10.509514426347591)};
        \addplot+[black, boxplot] coordinates {(5,-4.1568545135954) (5,-1.7498011555566764) (5,-0.1943435778486148) (5,-0.7964847051520418) (5,0.069721644427945) (5,-4.117744522123947) (5,-2.959464525504708) (5,0.0) (5,-0.8842393994836755) (5,0.5063792201799713)};
        \node[draw,fill=white,font=\scriptsize,anchor=north west,inner xsep=2pt,inner ysep=1pt] at (rel axis cs:0.03,0.97){\raisebox{0.5ex}{\textcolor{red}{\rule{0.2cm}{0.7pt}}}~F-S};
      \nextgroupplot[title={Work units (-)},title style={yshift=-4pt}, ymin=0.00, ymax=1.64, ytick distance=0.20]
      \draw[ red, thick] (axis cs:0.5,0.65) -- (axis cs:5.5,0.65);
        \addplot+[black, boxplot] coordinates {(1,0.1064204641560293) (1,0.1123280198252978) (1,0.1067210318536151) (1,0.1077748569816492) (1,0.1081343636646639) (1,0.1035256121157689) (1,0.1021565840819085) (1,0.1076152131975304) (1,0.1008475035819436) (1,0.104893698682955)};
        \addplot+[black, boxplot] coordinates {(2,0.1125350492262566) (2,0.1055298485597687) (2,0.1120909275803432) (2,0.1162684138961977) (2,0.1191874450288442) (2,0.1135183967365289) (2,0.1216576738269038) (2,0.1077632241258463) (2,0.1106387445925325) (2,0.1066682871234325)};
        \addplot+[black, boxplot] coordinates {(3,0.2569577204512302) (3,0.3668748556685651) (3,0.325121084660438) (3,0.2664353656557014) (3,0.276087011460943) (3,0.2491924996110069) (3,0.3466514509977199) (3,0.3156995139231398) (3,0.2646541242197003) (3,0.2879239006898461)};
        \addplot+[black, boxplot] coordinates {(4,1.4094580014546711) (4,0.5038486321767172) (4,0.47378706932067866) (4,0.48012638092041016) (4,0.46720170974731445) (4,0.4744281768798828) (4,1.0197128454844158) (4,1.0279446442921956) (4,0.4724737803141276) (4,1.0290738741556804)};
        \addplot+[black, boxplot] coordinates {(5,0.2594760103922955) (5,0.2310345967287875) (5,0.2283267680968511) (5,0.2428579681578494) (5,0.2314837711868793) (5,0.2168124867413813) (5,0.2707449659802469) (5,0.2220680543588323) (5,0.2280828257552806) (5,0.2238927822984589)};
        \node[draw,fill=white,font=\scriptsize,anchor=north west,inner xsep=2pt,inner ysep=1pt] at (rel axis cs:0.03,0.97){\raisebox{0.5ex}{\textcolor{red}{\rule{0.2cm}{0.7pt}}}~F-S};
    \end{groupplot}
  \end{tikzpicture}
  \captionsetup{skip=-3.5pt}
  \caption{Case study with 50h energy-to-power ratio storage.}
  \label{fig:ess50_solar}
\end{subfigure}
     \captionsetup{skip=-3.5pt}
        \caption{TSA performance relative to full-scale (F-S) optimization, with perturbed \textbf{PV} data in a cost minimization problem.}
        \label{fig:solar_comparison}
\end{figure}
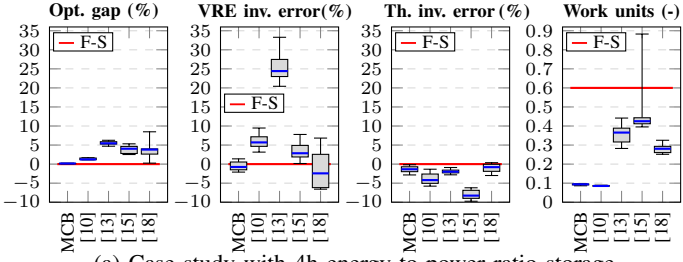
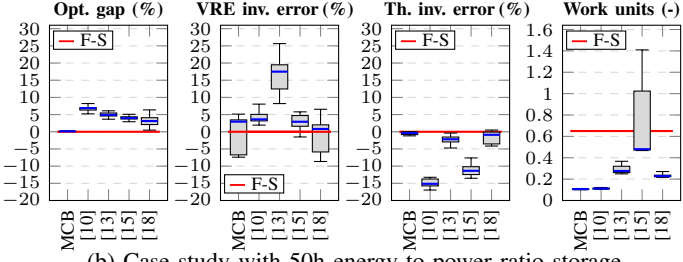

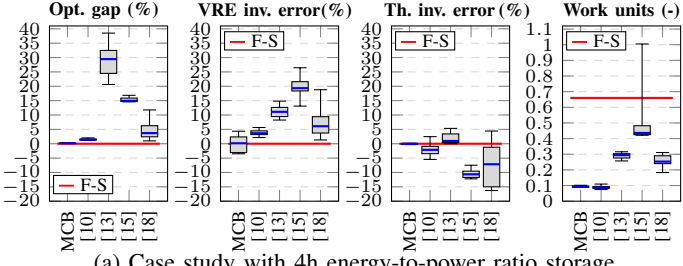
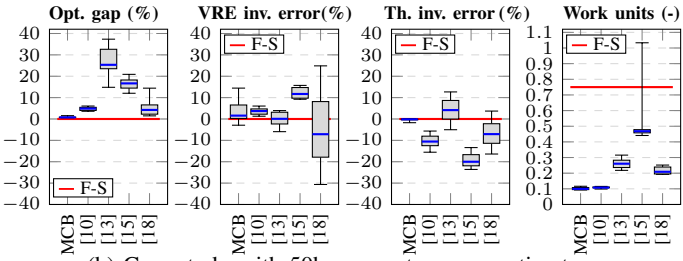
\begin{figure}[ht]

     \begin{subfigure}{\linewidth}
\hspace{-0.55cm}
  \centering
  \begin{tikzpicture}
    \begin{groupplot}[
      group style={
        group size=4 by 1,
        horizontal sep=0.7cm,
        x descriptions at=edge bottom,
      },
      width=0.355\linewidth,
      height=3.9cm,
      boxplot/draw direction=y,
      xtick={1,...,5},
      xticklabels={MCB,\cite{Pineda2018},\cite{Moradi2023},\cite{Zhang2023},\cite{Hilbers2023}},
      xticklabel style={rotate=90, anchor=east, font=\scriptsize},
      ymajorgrids,
      grid style={dashed, gray!30},
      boxplot/every median/.style={thick,draw=blue},      boxplot/whisker range=5000,
      title style={font=\scriptsize\bfseries, yshift=-3pt},
      tick label style={font=\scriptsize},
      boxplot/every box/.style={fill=gray!30},
      boxplot/draw/marker/.style={draw=none},
    ]
      \nextgroupplot[title={\!\!\!Opt. gap (\%)},title style={yshift=-4pt}, ymin=-20.00, ymax=41.00, ytick distance=5.00]
      \draw[red, thick] (axis cs:0.5,0) -- (axis cs:5.5,0);
        \addplot+[black, boxplot] coordinates {(1,0.1296383854473972) (1,0.1104058698658181) (1,0.1881145750727117) (1,0.23238409039626) (1,0.1671428711863161) (1,0.1941623386641795) (1,0.1322248107480022) (1,0.177285098365806) (1,0.1319453015796233) (1,0.161807425079745)};
        \addplot+[black, boxplot] coordinates {(2,1.8403816979396133) (2,1.4752534773880208) (2,1.5650688101422523) (2,1.6667646234943962) (2,1.310047844607373) (2,1.1277416112706595) (2,1.242749349450113) (2,1.4589149384730078) (2,1.2646502792678256) (2,2.039305413946652)};
        \addplot+[black, boxplot] coordinates {(3,34.65054891911748) (3,31.783322566509536) (3,20.63155573440636) (3,33.191123354878634) (3,29.45778644153581) (3,29.899407942438543) (3,38.506570741627165) (3,25.613691916082708) (3,23.36627103850028) (3,27.549433491803686)};
        \addplot+[black, boxplot] coordinates {(4,16.848769927839435) (4,14.692899932195184) (4,15.895510703422463) (4,16.615949248237705) (4,15.895120938683029) (4,14.822362186045382) (4,14.5771707264439) (4,14.69647533765032) (4,15.427640766308867) (4,14.943790999938985)};
        \addplot+[black, boxplot] coordinates {(5,11.767402288699692) (5,3.6929034473385287) (5,2.427948139475759) (5,7.437870759467813) (5,3.96224437458322) (5,0.9933294153285984) (5,3.073479163688593) (5,5.335305226271834) (5,2.4083238814510777) (5,7.205958385805706)};
        \node[draw,fill=white,font=\scriptsize,anchor=north west,inner xsep=2pt,inner ysep=1pt] at (rel axis cs:0.03,0.15){\raisebox{0.5ex}{\textcolor{red}{\rule{0.2cm}{0.7pt}}}~F-S};
      \nextgroupplot[title={\!\!\!VRE inv. error\!\! (\%)},title style={yshift=-4pt}, ymin=-20.00, ymax=41.00, ytick distance=5.00]
      \draw[red, thick] (axis cs:0.5,0) -- (axis cs:5.5,0);
        \addplot+[black, boxplot] coordinates {(1,-3.219594164883482) (1,0.9487920322587464) (1,3.5768674986156745) (1,2.740292166781173) (1,4.3571566632652665) (1,-3.497484485494721) (1,0.0285196197021952) (1,-2.955574142512619) (1,2.012762166760808) (1,0.1653724400149344)};
        \addplot+[black, boxplot] coordinates {(2,4.774955496389694) (2,2.9521026996603488) (2,3.6523360926934343) (2,5.366840652238186) (2,5.656275588395218) (2,3.4081264897400856) (2,3.5598033363255066) (2,4.186996435854117) (2,4.070337013400724) (2,2.159744780602332)};
        \addplot+[black, boxplot] coordinates {(3,8.762299446573786) (3,10.11597928549178) (3,12.249676301541474) (3,11.10949534906371) (3,10.182661375491238) (3,8.242437751225985) (3,11.920108963667907) (3,14.378713720134884) (3,14.839340339483616) (3,13.047605584389656)};
        \addplot+[black, boxplot] coordinates {(4,22.90097238185917) (4,18.94712775504061) (4,21.13838085168126) (4,26.448309320707317) (4,22.05361138608) (4,13.095776933366674) (4,20.66879013962745) (4,17.798956260058333) (4,19.343148854204404) (4,19.33919095788536)};
        \addplot+[black, boxplot] coordinates {(5,3.298766566283558) (5,4.535014206018634) (5,6.096823212253592) (5,4.155763533422778) (5,10.029652253565292) (5,6.441440853832287) (5,8.80287191528956) (5,18.800485470082208) (5,13.58295746131854) (5,1.3229722283829697)};
        \node[draw,fill=white,font=\scriptsize,anchor=north west,inner xsep=2pt,inner ysep=1pt] at (rel axis cs:0.03,0.97){\raisebox{0.5ex}{\textcolor{red}{\rule{0.2cm}{0.7pt}}}~F-S};
      \nextgroupplot[title={$\;$ Th. inv. error\! (\%)},title style={yshift=-4pt}, ymin=-20.00, ymax=41.00, ytick distance=5.00]
      \draw[red, thick] (axis cs:0.5,0) -- (axis cs:5.5,0);
        \addplot+[black, boxplot] coordinates {(1,-0.2974636613276565) (1,0.0) (1,-0.2587655054298893) (1,-0.1665077823575134) (1,0.0) (1,0.0890579575203954) (1,-0.0008579377820496) (1,0.1469027835495793) (1,-0.0653960448584761) (1,-0.01221828905733)};
        \addplot+[black, boxplot] coordinates {(2,-2.1244210177195235) (2,-0.9779926344844718) (2,-3.6342521260832648) (2,-5.467180764265264) (2,-0.8379976140814994) (2,-0.6814859533001643) (2,-2.175404854662711) (2,-1.6656539326189082) (2,2.467762669075834) (2,-3.286223482424667)};
        \addplot+[black, boxplot] coordinates {(3,0.797420459801412) (3,4.682539420161534e-13) (3,5.338599139275037) (3,3.390813992678708) (3,-5.543468005097635e-13) (3,4.425010281069922) (3,0.9687160187728596) (3,2.934474754929272) (3,0.772963922069202) (3,3.445208722368783)};
        \addplot+[black, boxplot] coordinates {(4,-11.360160445064595) (4,-12.217812951440685) (4,-7.476367564324914) (4,-11.065618889732812) (4,-10.642013960292822) (4,-10.07870199593546) (4,-9.369921213459053) (4,-7.885015079789009) (4,-12.111988045920066) (4,-9.708750260814652)};
        \addplot+[black, boxplot] coordinates {(5,-1.4974319736596895) (5,-14.309250670526838) (5,-10.171440493021937) (5,-16.28982592290948) (5,-6.339983259328394) (5,4.425010281069952) (5,-15.688412569004145) (5,-0.9344525004167202) (5,-0.441319749605906) (5,-7.134886352248156)};
        \node[draw,fill=white,font=\scriptsize,anchor=north west,inner xsep=2pt,inner ysep=1pt] at (rel axis cs:0.03,0.97){\raisebox{0.5ex}{\textcolor{red}{\rule{0.2cm}{0.7pt}}}~F-S};
      \nextgroupplot[title={Work units (-)},title style={yshift=-4pt}, ymin=0.00, ymax=1.12, ytick distance=0.10]
        \draw[ red, thick] (axis cs:0.5,0.66) -- (axis cs:5.5,0.66);
        \addplot+[black, boxplot] coordinates {(1,0.0911300754152541) (1,0.098052393511336) (1,0.0909126885913362) (1,0.0911669116042837) (1,0.0896385392638677) (1,0.1017368061375178) (1,0.1015244293480414) (1,0.0878615728475021) (1,0.096087969512999) (1,0.0900736992899965)};
        \addplot+[black, boxplot] coordinates {(2,0.10944520472782) (2,0.0841296791941395) (2,0.0948022478029033) (2,0.0745944647354257) (2,0.0817770155391714) (2,0.0923552859514226) (2,0.0866885131316774) (2,0.0897717605210022) (2,0.0968604754613079) (2,0.0825098495489768)};
        \addplot+[black, boxplot] coordinates {(3,0.2838106157341754) (3,0.3158266901041971) (3,0.3058466463835947) (3,0.2999712851311322) (3,0.2680284889578878) (3,0.2875313277153129) (3,0.3063285474612223) (3,0.2943292777265523) (3,0.3013816313388829) (3,0.2573918090610039)};
        \addplot+[black, boxplot] coordinates {(4,0.9719988505045573) (4,0.4277982711791992) (4,0.43630051612854) (4,1.0048339366912842) (4,0.4275425275166829) (4,0.4215552806854248) (4,0.5221668084462484) (4,0.4426408608754476) (4,0.44084469477335614) (4,0.42279283205668133)};
        \addplot+[black, boxplot] coordinates {(5,0.2529795553667607) (5,0.2485597903179423) (5,0.255502734898632) (5,0.233393369126376) (5,0.2791448741082199) (5,0.3087637778306408) (5,0.3111289115401618) (5,0.3011843543643787) (5,0.1829533196693579) (5,0.2466787637139277)};
        \node[draw,fill=white,font=\scriptsize,anchor=north west,inner xsep=2pt,inner ysep=1pt] at (rel axis cs:0.03,0.97){\raisebox{0.5ex}{\textcolor{red}{\rule{0.2cm}{0.7pt}}}~F-S};
    \end{groupplot}
  \end{tikzpicture}
  \captionsetup{skip=-3.5pt}
  \caption{Case study with 4h energy-to-power ratio storage.}
  \label{fig:ess4_wind}
\end{subfigure}
     \begin{subfigure}{\linewidth}
\hspace{-0.55cm}
  \centering
  \begin{tikzpicture}
    \begin{groupplot}[
      group style={
        group size=4 by 1,
        horizontal sep=0.7cm,
        x descriptions at=edge bottom,
      },
      width=0.355\linewidth,
      height=3.9cm,
      boxplot/draw direction=y,
      xtick={1,...,5},
      xticklabels={MCB,\cite{Pineda2018},\cite{Moradi2023},\cite{Zhang2023},\cite{Hilbers2023}},
      xticklabel style={rotate=90, anchor=east, font=\scriptsize},
      ymajorgrids,
      grid style={dashed, gray!30},
      boxplot/every median/.style={thick,draw=blue},      
      boxplot/whisker range=5000,
      title style={font=\scriptsize\bfseries, yshift=-3pt},
      tick label style={font=\scriptsize},
      boxplot/every box/.style={fill=gray!30},
      boxplot/draw/marker/.style={draw=none},
    ]
      \nextgroupplot[title={\!\!\!Opt. gap (\%)},title style={yshift=-4pt}, ymin=-40.00, ymax=42.00, ytick distance=10.00]
      \draw[red, thick] (axis cs:0.5,0) -- (axis cs:5.5,0);
        \addplot+[black, boxplot] coordinates {(1,0.7667534334671867) (1,1.5600888008757672) (1,0.7400517239650932) (1,0.7889435814690409) (1,0.9097266319047916) (1,0.6122589094122592) (1,0.7499589184854383) (1,0.7254369112151297) (1,0.7152944901126311) (1,0.8977864199380955)};
        \addplot+[black, boxplot] coordinates {(2,4.866953109548901) (2,4.116412548857692) (2,6.057057160107052) (2,4.541971169446871) (2,5.936218068321683) (2,4.930471558410782) (2,5.8576081856711175) (2,5.042548124273769) (2,3.612514199102633) (2,3.920219582870336)};
        \addplot+[black, boxplot] coordinates {(3,29.52491275386445) (3,25.344393451129783) (3,32.97504610703231) (3,24.85352852945514) (3,14.82693863644149) (3,37.34222291158808) (3,33.57258856283819) (3,25.04736875397112) (3,32.25484297789574) (3,22.185692964804534)};
        \addplot+[black, boxplot] coordinates {(4,12.034877771689274) (4,14.652756565842964) (4,20.8849910002664) (4,19.31366571567202) (4,17.172651907055123) (4,14.245453378449303) (4,14.962921702403008) (4,18.55967240619865) (4,17.93330647742903) (4,16.65510332461139)};
        \addplot+[black, boxplot] coordinates {(5,4.237799115782023) (5,4.960450144940239) (5,14.42344334786721) (5,1.462136429217885) (5,6.1727752748446125) (5,7.000278696112434) (5,4.064325272772976) (5,1.7753285180301848) (5,2.72646099281916) (5,14.352553593270438)};
        \node[draw,fill=white,font=\scriptsize,anchor=north west,inner xsep=2pt,inner ysep=1pt] at (rel axis cs:0.03,0.15){\raisebox{0.5ex}{\textcolor{red}{\rule{0.2cm}{0.7pt}}}~F-S};
      \nextgroupplot[title={\!\!\!VRE inv. error\!\! (\%)}, title style={yshift=-4pt},ymin=-40.00, ymax=42.00, ytick distance=10.00]
      \draw[red, thick] (axis cs:0.5,0) -- (axis cs:5.5,0);
        \addplot+[black, boxplot] coordinates {(1,2.3660164748063406) (1,14.446504709596756) (1,5.34054282848302) (1,-2.9123613414690084) (1,9.239218494030869) (1,1.0354207175927743) (1,1.5739347072916114) (1,1.088115840720578) (1,-0.53410202870775) (1,7.754830515262403)};
        \addplot+[black, boxplot] coordinates {(2,5.432724638975258) (2,2.1517038134720923) (2,3.7044714390622415) (2,3.7169122046627407) (2,1.289116852527901) (2,4.398587211725136) (2,6.053698391835019) (2,2.3325649289431274) (2,5.024018644823001) (2,3.52723330210067)};
        \addplot+[black, boxplot] coordinates {(3,-5.931008826070168) (3,3.9303342882890746) (3,0.0780845184000983) (3,3.1165230699508637) (3,-1.9972505740145765) (3,-0.6513873989612954) (3,3.197847065643464) (3,-2.670302727498428) (3,2.6821786470148665) (3,3.519375983353906)};
        \addplot+[black, boxplot] coordinates {(4,11.66192795565564) (4,9.929604605223364) (4,9.230605824013722) (4,9.522699452068895) (4,11.702901512001056) (4,14.929581086832542) (4,15.518438640027762) (4,11.76712374077651) (4,14.50938204690524) (4,15.711143705096694)};
        \addplot+[black, boxplot] coordinates {(5,-13.468870001940072) (5,-21.54939187032857) (5,-7.123432741873957) (5,8.390844304760755) (5,7.901849652285913) (5,13.186651111477708) (5,24.85929767416297) (5,2.4897355715268725) (5,-30.64247431290647) (5,-14.182542295211626)};

        \node[draw,fill=white,font=\scriptsize,anchor=north west,inner xsep=2pt,inner ysep=1pt] at (rel axis cs:0.03,0.97){\raisebox{0.5ex}{\textcolor{red}{\rule{0.2cm}{0.7pt}}}~F-S};
      \nextgroupplot[title={$\;$ Th. inv. error\! (\%)},title style={yshift=-4pt}, ymin=-40.00, ymax=42.00, ytick distance=10.00]
      \draw[red, thick] (axis cs:0.5,0) -- (axis cs:5.5,0);
        \addplot+[black, boxplot] coordinates {(1,-0.1614087898104217) (1,-1.737027135226758) (1,-0.370365914282899) (1,0.0973664775389483) (1,-0.4184465445572488) (1,-0.1194768205077546) (1,-0.0427393114373277) (1,-0.0375266128370417) (1,0.0483723489917836) (1,-0.2499531901140122)};
        \addplot+[black, boxplot] coordinates {(2,-10.558326375096303) (2,-7.213533260961637) (2,-12.648700205303625) (2,-8.862901180414388) (2,-15.547725312804852) (2,-10.245688394632092) (2,-12.322362636409936) (2,-11.595887975562473) (2,-5.6478523917571) (2,-6.7183899675248)};
        \addplot+[black, boxplot] coordinates {(3,7.69064239383655) (3,9.656150838685033) (3,2.2287820700344083) (3,12.666074584103129) (3,0.0553826169485122) (3,-5.02105287968814) (3,-0.0868357378703781) (3,4.155379211040576) (3,7.194361435108452) (3,10.021639190657508)};
        \addplot+[black, boxplot] coordinates {(4,-13.381313236292003) (4,-16.907778037750703) (4,-23.6033762962895) (4,-22.060779044708028) (4,-18.13612477282964) (4,-15.982108380911756) (4,-16.461165495369592) (4,-21.74968093238429) (4,-20.00808168964563) (4,-20.16544261549277)};
        \addplot+[black, boxplot] coordinates {(5,-4.428158797702399) (5,-8.213820982973283) (5,-14.099046832225614) (5,-8.661975298087157) (5,-6.394360539041903) (5,3.703359100601197) (5,-7.081374948141786) (5,-0.0037489186887781) (5,2.7752159358429216) (5,-16.392797828222)};

        \node[draw,fill=white,font=\scriptsize,anchor=north west,inner xsep=2pt,inner ysep=1pt] at (rel axis cs:0.03,0.97){\raisebox{0.5ex}{\textcolor{red}{\rule{0.2cm}{0.7pt}}}~F-S};
      \nextgroupplot[title={Work units (-)},title style={yshift=-4pt}, ymin=0.00, ymax=1.12, ytick distance=0.10]
      \draw[ red, thick] (axis cs:0.5,0.75) -- (axis cs:5.5,0.75);
        \addplot+[black, boxplot] coordinates {(1,0.0925351729754248) (1,0.1011568213850926) (1,0.0997293673075521) (1,0.1118792429536199) (1,0.09661976289496) (1,0.1067369027617644) (1,0.1099655288759106) (1,0.0973850223694426) (1,0.116369375510798) (1,0.0935140704833079)};
        \addplot+[black, boxplot] coordinates {(2,0.109286451921024) (2,0.1058050191167711) (2,0.1107413817127621) (2,0.098701802288294) (2,0.1127687739535558) (2,0.107269074399381) (2,0.1077350861651517) (2,0.1074051553663183) (2,0.1154418071599939) (2,0.1134757974893808)};
        \addplot+[black, boxplot] coordinates {(3,0.2912729501560712) (3,0.2408505655160465) (3,0.2702961085920465) (3,0.2475729376830079) (3,0.2890326389572943) (3,0.2605753770131477) (3,0.2765975558183468) (3,0.3154376268986351) (3,0.2180851260850629) (3,0.231282505758855)};
        \addplot+[black, boxplot] coordinates {(4,0.44086678822835285) (4,0.5312932332356771) (4,0.474340279897054) (4,1.0333775679270427) (4,0.4819819927215576) (4,0.4682099024454753) (4,0.4713714917500813) (4,0.45578885078430176) (4,0.45968198776245117) (4,0.46747581164042157)};
        \addplot+[black, boxplot] coordinates {(5,0.2088587966420859) (5,0.2405097808344113) (5,0.1961500207194951) (5,0.2518481491208304) (5,0.1926528128268449) (5,0.2467877518552105) (5,0.2381451066651806) (5,0.1985667987818286) (5,0.2260111087110588) (5,0.1939878307748851)};
        \node[draw,fill=white,font=\scriptsize,anchor=north west,inner xsep=2pt,inner ysep=1pt] at (rel axis cs:0.03,0.97){\raisebox{0.5ex}{\textcolor{red}{\rule{0.2cm}{0.7pt}}}~F-S};
    \end{groupplot}
  \end{tikzpicture}
  \captionsetup{skip=-3.5pt}
  \caption{Case study with 50h energy-to-power ratio storage.}
  \label{fig:ess50_wind}
\end{subfigure}
     \captionsetup{skip=-3.5pt}
        \caption{TSA performance relative to full-scale (F-S) optimization, with perturbed \textbf{wind} data in a cost minimization problem.}
        \label{fig:wind_comparison}
\end{figure}

Figures \ref{fig:solar_comparison} and \ref{fig:wind_comparison} show TSA method sensitivity across 10 scenarios with uniform data perturbations. 
Models include either a PV (Fig.~\ref{fig:solar_comparison}) or a wind unit (Fig.~\ref{fig:wind_comparison}), each coupled with storage ($\overline{E}_s / \overline{P}^\mathrm{d}_s$ ratios of 4--50 h). 
Performance is evaluated model via optimality gap, solution time (work units), and VRE/thermal investment errors relative to full-scale model.

Our MCB TSA method outperforms benchmarks, reducing solution time tenfold compared to the full-scale optimization. 
By leveraging structural information from the full-scale model, it demonstrates superior robustness to data perturbations compared to input-dependent alternatives.
Although validated on a small case study, the method generalizes to larger systems provided that marginal costs are accurately captured. Future work will focus on aggregating heterogeneous marginal costs in systems with diverse technologies and network constraints.

\subsection{Market Participation Problem}\label{sub:MP}

\begin{figure}[ht]
    \begin{subfigure}{\linewidth}
\hspace{-0.55cm}
  \centering
  \begin{tikzpicture}
    \begin{groupplot}[
      group style={
        group size=4 by 1,
        horizontal sep=0.7cm,
        x descriptions at=edge bottom,
      },
      width=0.355\linewidth,
      height=3.9cm,
      boxplot/draw direction=y,
      xtick={1,...,5},
      xticklabels={MCB,\cite{Pineda2018},\cite{Moradi2023},\cite{Zhang2023},\cite{Hilbers2023}},
      xticklabel style={rotate=90, anchor=east, font=\scriptsize},
      ymajorgrids,
      grid style={dashed, gray!30},
      boxplot/every median/.style={thick,draw=blue},      
      boxplot/whisker range=5000,
      title style={font=\scriptsize\bfseries, yshift=-3pt},
      tick label style={font=\scriptsize},
      boxplot/every box/.style={fill=gray!30},
      boxplot/draw/marker/.style={draw=none},
    ]
      \nextgroupplot[title={Opt. Gap(\%)},title style={yshift=-4pt}, ymin=-6.00, ymax=14.40, ytick distance=2.00]
      \draw[red, thick] (axis cs:0.5,0) -- (axis cs:5.5,0);
        \addplot+[black, boxplot] coordinates {(1,1.5477589413239172) (1,2.351217961882106) (1,1.2803272353098576) (1,0.3437126920792455) (1,0.5099968539202829) (1,0.6521966103516992) (1,0.2837334879286612) (1,0.6188716266187689) (1,0.194759807476902) (1,0.3005178286579989)};
        \addplot+[black, boxplot] coordinates {(2,0.5394442984711758) (2,0.7103694173542178) (2,0.6618140316868233) (2,0.6630125091060366) (2,0.6099623286961163) (2,0.3765974134574911) (2,0.4399193017940381) (2,0.788591359989832) (2,0.4350105314912234) (2,0.9209348227829692)};
        \addplot+[black, boxplot] coordinates {(3,3.7834215528992825) (3,3.4833920413715034) (3,3.4311040100660843) (3,3.297407828936189) (3,2.884285781787685) (3,3.0455087496038518) (3,3.230346142326016) (3,3.136005090463362) (3,4.377239481287198) (3,3.566113589819299)};
        \addplot+[black, boxplot] coordinates {(4,1.3506057466341566) (4,2.6570192293259622) (4,3.972538866910178) (4,1.694859359084594) (4,1.6835232484356446) (4,1.228680178272253) (4,2.244982423102914) (4,1.5625463549570209) (4,1.0786601113938707) (4,1.67670173509482)};
        \addplot+[black, boxplot] coordinates {(5,1.7740566771482489) (5,5.044857797614932) (5,0.8947622190338842) (5,0.3644353802943403) (5,7.639467184777445) (5,2.799449097948669) (5,12.363772654673015) (5,8.448817290831455) (5,9.925902019348094) (5,9.39744797196027)};
        \node[draw,fill=white,font=\scriptsize,anchor=north west,inner xsep=2pt,inner ysep=1pt] at (rel axis cs:0.03,0.97){\raisebox{0.5ex}{\textcolor{red}{\rule{0.2cm}{0.7pt}}}~F-S};
      \nextgroupplot[title={\!\!\!VRE inv. error \!\!(\%)}, title style={yshift=-4pt},ymin=-6.00, ymax=14.40, ytick distance=2.00]
      \draw[red, thick] (axis cs:0.5,0) -- (axis cs:5.5,0);
        \addplot+[black, boxplot] coordinates {(1,-2.591651598708497) (1,-2.5369503861104645) (1,-2.939944593725099) (1,-2.0293518280886387) (1,-3.102007373626629) (1,-3.0770096041928623) (1,-2.1512728674952406) (1,-2.8064515985360674) (1,-2.6076295347800453) (1,-2.843849302258782)};
        \addplot+[black, boxplot] coordinates {(2,-2.041856407596092) (2,0.0900350021473942) (2,-0.2531630239426036) (2,0.6133402972953511) (2,-0.8540011179956782) (2,-2.523474258180015) (2,0.2700990336941851) (2,-0.8080812961985775) (2,-1.469661818985054) (2,-0.3685739323548238)};
        \addplot+[black, boxplot] coordinates {(3,-2.863072857047021) (3,-2.235190566467478) (3,-2.604784648452412) (3,-2.016843210666648) (3,-3.524804769492687) (3,-3.112702207380062) (3,-2.029560554000576) (3,-2.8777465930531627) (3,-3.1094818640498274) (3,-1.733194628508374)};
        \addplot+[black, boxplot] coordinates {(4,-0.5158439097810887) (4,1.2890760237582426) (4,1.6114622846366495) (4,1.0636783767970164) (4,0.4773542466858992) (4,0.3073884029172979) (4,1.1598465170814312) (4,0.6468075933906545) (4,0.7570343240882609) (4,0.7762129786292083)};
        \addplot+[black, boxplot] coordinates {(5,-2.863072857047021) (5,-2.5216065986998286) (5,-2.804454181781332) (5,-2.0293518280886387) (5,-3.223493978608856) (5,-3.5189584591609457) (5,-1.8334330837798356) (5,-2.5370952631250057) (5,-2.966617028211472) (5,-1.961844511873585)};
        \node[draw,fill=white,font=\scriptsize,anchor=north west,inner xsep=2pt,inner ysep=1pt] at (rel axis cs:0.03,0.97){\raisebox{0.5ex}{\textcolor{red}{\rule{0.2cm}{0.7pt}}}~F-S};
      \nextgroupplot[title={$\;$ Th. inv. error\! (\%)},title style={yshift=-4pt}, ymin=-6.00, ymax=14.40, ytick distance=2.00]
      \draw[red, thick] (axis cs:0.5,0) -- (axis cs:5.5,0);
        \addplot+[black, boxplot] coordinates {(1,7.972744614142127) (1,7.804466286011413) (1,9.044204644319768) (1,6.242931675560126) (1,9.542761300722503) (1,9.46586021119414) (1,6.617999570783377) (1,8.633537732552234) (1,8.021897827414797) (1,8.748584963856464)};
        \addplot+[black, boxplot] coordinates {(2,6.281399739311307) (2,-0.2769763029925299) (2,0.7788099822692489) (2,-1.8868298325032853) (2,2.627179061168732) (2,7.763009430301936) (2,-0.8309105349050293) (2,2.4859150841365727) (2,4.521147193458567) (2,1.1338506439520846)};
        \addplot+[black, boxplot] coordinates {(3,8.80772273259778) (3,6.876157103549521) (3,8.013146052231) (3,6.20445118989941) (3,10.843420564662246) (3,9.575662010929577) (3,6.243573782852234) (3,8.852863811657445) (3,9.565755210588588) (3,5.331857934371398)};
        \addplot+[black, boxplot] coordinates {(4,1.5868999349660386) (4,-3.965607850515149) (4,-4.957370526629662) (4,-3.272212998851413) (4,-1.4684934893253414) (4,-0.9456244948317996) (4,-3.5680567854494027) (4,-1.989785879846472) (4,-2.328878362007981) (4,-2.3878780033077107)};
        \addplot+[black, boxplot] coordinates {(5,8.80772273259778) (5,7.75726391571611) (5,8.627393043319517) (5,6.242931675560126) (5,9.916492737480874) (5,10.8254354546138) (5,5.640223304467376) (5,7.80491197378754) (5,9.126257536191194) (5,6.035257699614827)};
        \node[draw,fill=white,font=\scriptsize,anchor=north west,inner xsep=2pt,inner ysep=1pt] at (rel axis cs:0.03,0.97){\raisebox{0.5ex}{\textcolor{red}{\rule{0.2cm}{0.7pt}}}~F-S};
      \nextgroupplot[title={Work units (-)},title style={yshift=-4pt}, ymin=0.00, ymax=1.64, ytick distance=0.20]
      \draw[ red, thick] (axis cs:0.5,0.46) -- (axis cs:5.5,0.46);
        \addplot+[black, boxplot] coordinates {(1,0.0545220139198887) (1,0.0736502328016992) (1,0.0990142025533757) (1,0.0520153921046526) (1,0.0724731868789125) (1,0.1086807855024936) (1,0.0557509673611977) (1,0.0619226882399836) (1,0.0584591416870346) (1,0.0627013660635562)};
        \addplot+[black, boxplot] coordinates {(2,0.0688641691687747) (2,0.0606440955903246) (2,0.1112172723044775) (2,0.0589886959387814) (2,0.0900619738219292) (2,0.0670586447691558) (2,0.0583027379623838) (2,0.0586653517602891) (2,0.0589648165802114) (2,0.0461956394384728)};
        \addplot+[black, boxplot] coordinates {(3,0.3792998837237205) (3,0.4099131168820756) (3,0.4777800683659266) (3,0.5871140592783985) (3,0.5643167936387041) (3,0.4074046072016569) (3,0.4024027885809679) (3,0.3724031206856235) (3,0.4091751699486257) (3,0.4174368801034228)};
        \addplot+[black, boxplot] coordinates {(4,0.12713805805731176) (4,0.19375945731270802) (4,0.11720031854961743) (4,0.7292420993001519) (4,0.10985619233614753) (4,0.11502782412549456) (4,0.18432706283264141) (4,0.11393931531535094) (4,0.1632479738352486) (4,0.1455045670592002)};
        \addplot+[black, boxplot] coordinates {(5,0.4067308305666695) (5,0.3617834563380108) (5,0.38483240464231) (5,1.51412885798554) (5,1.5428607578853033) (5,0.4263959236887064) (5,0.3435111899435482) (5,0.3871085735736201) (5,0.3133008004501366) (5,0.4222795244788299)};

        \node[draw,fill=white,font=\scriptsize,anchor=north west,inner xsep=2pt,inner ysep=1pt] at (rel axis cs:0.03,0.97){\raisebox{0.5ex}{\textcolor{red}{\rule{0.2cm}{0.7pt}}}~F-S};
    \end{groupplot}
  \end{tikzpicture}
  \captionsetup{skip=-3.5pt}
  \caption{Case study with PV as the sole renewable energy source.}
  \label{fig:ess4_solar_mp}
\end{subfigure}
    \begin{subfigure}{\linewidth}
\hspace{-0.55cm}
  \centering
  \begin{tikzpicture}
    \begin{groupplot}[
      group style={
        group size=4 by 1,
        horizontal sep=0.7cm,
        x descriptions at=edge bottom,
      },
      width=0.355\linewidth,
      height=3.9cm,
      boxplot/draw direction=y,
      xtick={1,...,5},
      xticklabels={MCB,\cite{Pineda2018},\cite{Moradi2023},\cite{Zhang2023},\cite{Hilbers2023}},
      xticklabel style={rotate=90, anchor=east, font=\scriptsize},
      ymajorgrids,
      grid style={dashed, gray!30},
      boxplot/every median/.style={thick,draw=blue},      boxplot/whisker range=5000,
      title style={font=\scriptsize\bfseries, yshift=-3pt},
      tick label style={font=\scriptsize},
      boxplot/every box/.style={fill=gray!30},
      boxplot/draw/marker/.style={draw=none},
    ]
      \nextgroupplot[title={Opt. gap (\%)},title style={yshift=-4pt}, ymin=-20.00, ymax=102.00, ytick distance=10.00]
      \draw[red, thick] (axis cs:0.5,0) -- (axis cs:5.5,0);
        \addplot+[black, boxplot] coordinates {(1,13.283088257756065) (1,8.975328165588577) (1,8.450965407169816) (1,8.365797479367224) (1,8.614925780238467) (1,7.039771362204809) (1,8.146382436815685) (1,8.693054270686478) (1,9.204891410766146) (1,10.193803383420112)};
        \addplot+[black, boxplot] coordinates {(2,3.9411073714694433) (2,1.8818173808220655) (2,2.913719267245683) (2,2.663416099939937) (2,2.1836261241375508) (2,3.615303757118863) (2,5.18136212192013) (2,2.081760247189054) (2,5.7183732092468365) (2,0.5257684730231487)};
        \addplot+[black, boxplot] coordinates {(3,21.763605761124023) (3,9.413847799874988) (3,16.78403316981537) (3,16.485541055042244) (3,5.827003698465474) (3,40.46256590436646) (3,27.507402033434435) (3,38.84597110499857) (3,13.03582260800334) (3,5.6034039561569084)};
        \addplot+[black, boxplot] coordinates {(4,35.20632379503792) (4,35.853563089222696) (4,25.602511955333835) (4,40.72880423179755) (4,27.287340948202765) (4,28.029154537610054) (4,23.872291430805294) (4,37.41944684418156) (4,32.52443897272157) (4,36.34016258820652)};
        \addplot+[black, boxplot] coordinates {(5,64.32909345546113) (5,29.12938471187873) (5,30.329498171032665) (5,11.59540021506736) (5,20.3884240497864) (5,92.52962547850112) (5,33.08119712701314) (5,28.269871878288004) (5,18.791260154507423) (5,42.14370150310162)};
        \draw[red, thick] (axis cs:0.5,0) -- (axis cs:5.5,0);
        \node[draw,fill=white,font=\scriptsize,anchor=north west,inner xsep=2pt,inner ysep=1pt] at (rel axis cs:0.03,0.97){\raisebox{0.5ex}{\textcolor{red}{\rule{0.2cm}{0.7pt}}}~F-S};
      \nextgroupplot[title={\!\!\!VRE inv. error\!\! (\%)},title style={yshift=-4pt}, ymin=-20.00, ymax=102.00, ytick distance=10.00]
      \draw[red, thick] (axis cs:0.5,0) -- (axis cs:5.5,0);
        \addplot+[black, boxplot] coordinates {(1,2.4541935375604638) (1,0.5931388665716232) (1,1.0481976632522914) (1,-0.1833446976195358) (1,0.9252468533295162) (1,0.487320902947452) (1,-0.4141275229966507) (1,1.278196698856658) (1,2.38449705137545) (1,0.5310956851927312)};
        \addplot+[black, boxplot] coordinates {(2,-0.3333992875936151) (2,0.6920503488418782) (2,0.4579063413581461) (2,0.6914548596941529) (2,0.947563656899136) (2,1.122463942247025) (2,-0.2154148773326173) (2,1.324026516522548) (2,0.7079007723694639) (2,0.5347278253843893)};
        \addplot+[black, boxplot] coordinates {(3,-3.263065139550498) (3,-2.5531706539797665) (3,-1.790950049178836) (3,-2.765735279967248) (3,-0.6270136797004741) (3,0.8104008467992323) (3,-2.5446059025123144) (3,0.8092569528197143) (3,-1.9797894450286375) (3,-0.5406587875734251)};
        \addplot+[black, boxplot] coordinates {(4,4.75858363876938) (4,5.299774914138779) (4,1.918654093821108) (4,6.175987952244598) (4,1.5111097439810517) (4,2.6683818024117) (4,1.580771703299716) (4,5.240761186013051) (4,4.847306859171907) (4,4.658293679536965)};
        \addplot+[black, boxplot] coordinates {(5,4.789611605291517) (5,4.219811257067732) (5,1.8203595506817132) (5,-1.407902987227016) (5,-0.777837710266228) (5,-0.0831236190644485) (5,-1.302941186101393) (5,0.4087993836322238) (5,-2.38213284823676) (5,-0.8144606359678856)};
        \node[draw,fill=white,font=\scriptsize,anchor=north west,inner xsep=2pt,inner ysep=1pt] at (rel axis cs:0.03,0.97){\raisebox{0.5ex}{\textcolor{red}{\rule{0.2cm}{0.7pt}}}~F-S};
      \nextgroupplot[title={$\;$ Th. inv. error\! (\%)},title style={yshift=-4pt}, ymin=-20.00, ymax=102.00, ytick distance=10.00]
      \draw[red, thick] (axis cs:0.5,0) -- (axis cs:5.5,0);
        \addplot+[black, boxplot] coordinates {(1,-7.765298386126837) (1,-1.8767469691555376) (1,-3.3165956548340656) (1,0.580119808295649) (1,-2.9275677679725582) (1,-1.5419290138565351) (1,1.3103382992252937) (1,-4.044334161457612) (1,-7.544772171136802) (1,-1.6804365279216829)};
        \addplot+[black, boxplot] coordinates {(2,1.0549065956957837) (2,-2.189712170101388) (2,-1.4488585839404384) (2,-2.1878279866224) (2,-2.998180226236497) (2,-3.551581122603969) (2,0.6815928628683031) (2,-4.189344000213459) (2,-2.23986439581405) (2,-1.6919289599311362)};
        \addplot+[black, boxplot] coordinates {(3,10.32464395092718) (3,8.078471259672236) (3,5.666733822608929) (3,8.751045660128419) (3,1.983930053005308) (3,-2.564184238713838) (3,8.051371583249054) (3,-2.560564850944587) (3,6.264239371127083) (3,1.7106950802105545)};
        \addplot+[black, boxplot] coordinates {(4,-15.056604658456012) (4,-16.76898457996346) (4,-6.070801389646969) (4,-19.54140476059427) (4,-4.781292867334638) (4,-8.443010132132788) (4,-5.001709836083269) (4,-16.582259612775438) (4,-15.337333243899453) (4,-14.739277827196236)};
        \addplot+[black, boxplot] coordinates {(5,-15.154780052805156) (5,-13.3518783432379) (5,-5.759788241937318) (5,4.454737015323024) (5,2.461151422557744) (5,0.2630109219551901) (5,4.122628038418155) (5,-1.293479566865945) (5,7.537291610807578) (5,2.5770297921703658)};
        \node[draw,fill=white,font=\scriptsize,anchor=north west,inner xsep=2pt,inner ysep=1pt] at (rel axis cs:0.03,0.97){\raisebox{0.5ex}{\textcolor{red}{\rule{0.2cm}{0.7pt}}}~F-S};
      \nextgroupplot[title={Work units (-)},title style={yshift=-4pt}, ymin=0.00, ymax=0.62, ytick distance=0.10]
      \draw[ red, thick] (axis cs:0.5,0.55) -- (axis cs:5.5,0.55);
        \addplot+[black, boxplot] coordinates {(1,0.0597981343406026) (1,0.0805577564990109) (1,0.0617280687912328) (1,0.0589698354309405) (1,0.0721951741445824) (1,0.0757503426014022) (1,0.0648320175093997) (1,0.0720037238194251) (1,0.0892137013428516) (1,0.0736816633454218)};
        \addplot+[black, boxplot] coordinates {(2,0.0809842013818065) (2,0.0560042846460726) (2,0.1082908601233401) (2,0.0944575751705229) (2,0.0832834384641313) (2,0.074581682405256) (2,0.0928826035051594) (2,0.0733406172445289) (2,0.0812148912572025) (2,0.1017353539382933)};
        \addplot+[black, boxplot] coordinates {(3,0.349642675033382) (3,0.3365800020763335) (3,0.3169415604823575) (3,0.3762061355858392) (3,0.3319715713879489) (3,0.3413097991211071) (3,0.4007503416033563) (3,0.3065706227161781) (3,0.3330013787680979) (3,0.3858877240556013)};
        \addplot+[black, boxplot] coordinates {(4,0.09422922461216356) (4,0.10707367769923397) (4,0.11062506378339793) (4,0.10615621840632423) (4,0.10531927255671843) (4,0.10147164418630904) (4,0.3884778383705787) (4,0.12268927710403497) (4,0.11388855360544374) (4,0.09701916325162609)};
        \addplot+[black, boxplot] coordinates {(5,0.4500184063723815) (5,0.3797493989986946) (5,0.3941982062306036) (5,0.5013553597447012) (5,0.4154804033843329) (5,0.4712698637419232) (5,0.5330329345995584) (5,0.4038973441520049) (5,0.377316115520255) (5,0.5893529090454055)};
        \node[draw,fill=white,font=\scriptsize,anchor=north west,inner xsep=2pt,inner ysep=1pt] at (rel axis cs:0.03,0.80){\raisebox{0.5ex}{\textcolor{red}{\rule{0.2cm}{0.7pt}}}~F-S};
    \end{groupplot}
  \end{tikzpicture}
  \captionsetup{skip=-3.5pt}
  \caption{Case study with wind as the sole renewable energy source.}
  \label{fig:ess4_wind_mp}
\end{subfigure}
    \captionsetup{skip=-3.5pt}
    \caption{TSA performance relative to full-scale (F-S) optimization, under perturbed \textbf{PV} and \textbf{wind} data for GEP with market participation and a storage unit with 4h energy-to-power ratio.}
\label{fig:solar_wind_comparison_mp}
\end{figure}
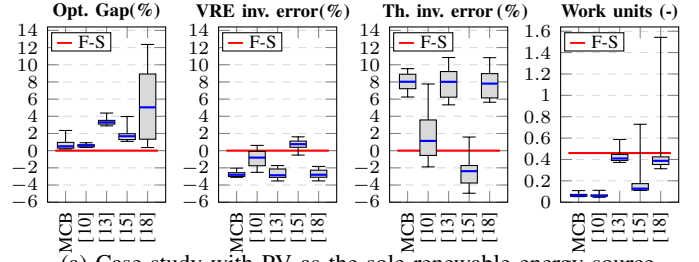
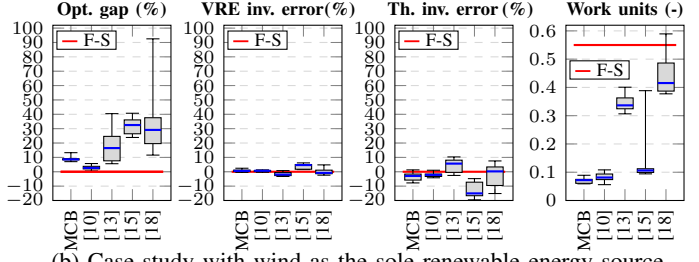

This subsection presents the performance of our TSA method for a GEP model incorporating market participation.

Fig. \ref{fig:solar_wind_comparison_mp} reports the performance of the analyzed TSA methods when both wind and PV capacity factor time series are perturbed with the aforementioned uniform noise. 
The inclusion of market prices degrades the performance of our MCB TSA method (see the optimality gaps in Fig.~\ref{fig:solar_wind_comparison_mp}) relative to the cases illustrated in the previous subsection, as the number of distinct marginal cost values in the GEP model increases.

Our MCB TSA method maintains the lowest optimality gap, though all methods yield similar investment errors. 
This suggests our approach excels at capturing operational decisions, as true marginal costs are more accurately reflecting operations rather than investments.
Consistent with previous results, our method reduces solution time tenfold compared to full-scale optimization.

\section{Conclusions}\label{s:con}
This paper investigates the GEP problem for a power system comprising renewable, thermal, and storage technologies, jointly optimizing market participation as well as operational and investment costs.
To overcome the heuristic nature of traditional TSA procedures, we propose a novel iterative TSA method that enables the explicit assessment of the optimality gap at each iteration.
Notably, the method targets exact temporal aggregation by employing ML-based estimates of the marginal costs as clustering features to construct an aggregated model that retains the same active constraints as its full-scale counterpart.
Numerical results show that the proposed method outperforms existing TSA methods in solution accuracy while achieving an order-of-magnitude runtime reduction relative to traditional full-scale optimization.
Future work will consider large-scale GEP problems and incorporate additional technical constraints, such as ramping and network constraints, to evaluate the applicability of the proposed MCB TSA method to a broader class of GEP models.

\appendices
\section{Full-scale Model}\label{sub:Full_scale}
The aggregated model presented in Subsection \ref{sub:Agg_model} collapses to the full-scale model when each representative time step corresponds to exactly one original time step (i.e., $\mathcal{R} = \mathcal{T}$).
The solution to the full-scale GEP model serves as benchmark for the true optimum.

The full-scale decision set is defined as:
\begin{equation} \label{eqn:decset}
\boldsymbol{z} \coloneqq \left\{x_g, p_{g,t}, p^{\mathrm{d}}_{s,t}, p^\mathrm{c}_{s,t}, e^\mathrm{ns}_t, e_{s,t}, o_t\right\}_{g\in\mathcal{G},t\in\mathcal{T},s\in\mathcal{S}},
\end{equation}
with the objective function:
\begin{align} \label{eqn:obj}
\!\!\!\!\!f\left(\boldsymbol{z}\right) \!\coloneqq\! \sum_{g \in \mathcal{G}} C_g^\mathrm{inv} x_g + \sum_{t \in \mathcal{T}} \bigg(& \! \Delta \sum_{g \in \mathcal{G}} C_g p_{g,t} + \Delta  \sum_{s \in \mathcal{S}} C_s^\mathrm{d} p_{s,t}^\mathrm{d}\nonumber \\ 
&+ C^\mathrm{ns} e_t^\mathrm{ns} - \pi_t o_t \bigg).
\end{align} 
 The \textbf{full-scale GEP model} is formulated as follows:

\begin{subequations}\label{eqn:model}
\begin{align}
    \!\!\min_{\boldsymbol{z} \geq \boldsymbol{0}} \; & f\left(\boldsymbol{z}\right)  \\
    \text{s.t.} \; &\sum_{g \in \mathcal{G}} p_{g,t} \Delta + \!\sum_{s \in \mathcal{S}} \left(p_{s,t}^\mathrm{d} \!-\! p_{s,t}^\mathrm{c}\right) \Delta \!+\! e_t^\mathrm{ns}\! - o_t \! =\! D_t, \; \forall t \label{eqn:balance} \\
     & p_{g,t} \leq F_{g,t} \, x_g, \; \forall g, \forall t, \label{eqn:gen_limits} \\
     & e_{s,t} = e_{s,t-1} + \left(\eta_s^\mathrm{c} p_{s,t}^\mathrm{c} - \frac{p_{s,t}^\mathrm{d}}{\eta_s^\mathrm{d}}\right) \Delta, \; \forall s ,\forall t, \label{eqn:soc_dynamic} \\
     & e_{s,0} = e_{s,T} = \underline{E}_s, \; \forall s, \label{eqn:soc_end} \\
     & \underline{E}_s \leq e_{s,t} \leq \overline{E}_s, \; \forall s, \forall t, \label{eqn:soc_bounds} \\
     & p_{s,t}^\mathrm{c} \leq \overline{P}_s^\mathrm{c},\; \forall s, \forall t, \label{eqn:storage_charging} \\ 
     & p_{s,t}^\mathrm{d} \leq \overline{P}_s^\mathrm{d}, \; \forall s, \forall t, \label{eqn:storage_discharging} \\
     & \sum_{g \in \mathcal{G}} C_g^\mathrm{inv} x_g \leq \overline{B}. \label{eqn:inv_limits}
\end{align}
\end{subequations}

In \eqref{eqn:model}, the constraints \eqref{eqn:balance} ensure the energy balance,
while the constraints \eqref{eqn:gen_limits} enforce the generation capacity limits.
The storage dynamics are described by the constraints \eqref{eqn:soc_dynamic}--\eqref{eqn:soc_bounds}.
Finally, the constraints \eqref{eqn:storage_charging}, \eqref{eqn:storage_discharging}, and \eqref{eqn:inv_limits} impose the charging, discharging, and investment limits, respectively.

\end{document}